\input amstex
\input epsf.tex
\documentstyle{amsppt}\nologo\footline={}

\def\PU{\mathop{\text{\rm PU}}}
\def\Aut{\mathop{\text{\rm Aut}}}
\def\T{\mathop{\text{\rm T}}}
\def\B{\mathop{\text{\rm B}}}
\def\S{\mathop{\text{\rm S}}}
\def\G{\mathop{\text{\rm G}}}
\def\Area{\mathop{\text{\rm Area}}}
\def\SU{\mathop{\text{\rm SU}}}
\def\Im{\mathop{\text{\rm Im}}}
\def\dist{\mathop{\text{\rm dist}}}

\hsize450pt

\topmatter\title A Hyperelliptic View on Teichm\"uller Space.
I\endtitle\author Sasha Anan$'$in and Eduardo C.~Bento
Gon\c calves\endauthor\thanks Second author supported by FAPESP,
project No.~04/15521\endthanks\address Departamento de Matem\'atica,
IMECC, Universidade Estadual de
Campinas,\newline13083-970--Campinas--SP, Brasil\endaddress\email
Ananin$_-$Sasha\@yahoo.com\endemail\address Departamento de
Matem\'atica, IMECC, Universidade Estadual de
Campinas,\newline13083-970--Campinas--SP, Brasil\endaddress\email
EduardoCBG\@gmail.com\endemail\subjclass 30F60 (51M10,
57S30)\endsubjclass\abstract We explicitly describe the Teichm\"uller
space $\Cal TH_n$ of hyperelliptic surfaces in terms of natural and
effective coordinates as the space of certain $(2n-6)$-tuples of
distinct points on the ideal boundary of the Poincar\'e disc. We
essentially use the concept of a simple earthquake which is a
particular case of a Fenchel-Nielsen twist deformation. Such
earthquakes generate a group that acts transitively on $\Cal TH_n$.
This fact can be interpreted as a continuous analog of the well-known
Dehn theorem saying that the mapping class group is generated by Dehn
twists. We find a simple and effective criterion that verifies if a
given representation of the surface group $\pi_1\Sigma$ in the group of
isometries of the hyperbolic plane is faithful and discrete.
The~article also contains simple and elementary proofs of several known
results, for instance, of~W.~M.~Goldman's theorem [Gol1] characterizing
the faithful discrete representations as having maximal Toledo
invariant (which is essentially the area of the representation in the
two-dimensional case).\endabstract\endtopmatter\document

\centerline{\bf1.~Introduction}

\medskip

This article is an attempt to an elementary study of Teichm\"uller
spaces and we hope it does not require from the reader any specific
knowledge in the field. We try to avoid the analytic methods typical in
the classic theory and worry more about the way of the proofs than
about the facts {\it per se,} having no prejudice against proving
well-known ones. Such elementary approach is motivated by its possible
extension to complex hyperbolic Teichm\"uller spaces and originates
from [Ana1].

Let $\Sigma=\Bbb D/\pi_1\Sigma$ be a hyperelliptic Riemann surface of
genus $g\ge2$, where $\Bbb D$ stands for the Poincar\'e disc. It is
well known [Mac] (and proven in Proposition 4.1) that the extension
$H_n$ of the fundamental group $\pi_1\Sigma$ with an isometry of
$\Bbb D$ induced by the hyperelliptic involution of $\Sigma$ is a group
with generators $r_1,\dots,r_n$ and defining relations
$r_n\dots r_1=1$, $r_i^2=1$, where $n=2g+2$. Moreover, every $r_i$ is a
reflection in some point $q_i\in\Bbb D$. In other words, a
hyperelliptic surface can be described as a certain geometric
configuration of $n$ points.

The following two concepts are crucial in this article. As is easy to
see, while moving the points $q_{i-1}$ and $q_i$ along the geodesic
they generate and preserving the distance between these two points,
new~configurations provide new hyperelliptic surfaces, i.e., the
relation $r_n\dots r_1=1$ remains valid. We call such a deformation a
{\it simple earthquake\/} (SE for short). This concept is nothing more
than a particular case of a Fenchel-Nielsen twist deformation [ImT]. It
appears naturally in the context of [Ana1]. The earthquake group
$\Cal E_n$, i.e., the formal group generated by the SEs, acts on the
Teichm\"uller space $\Cal TH_n$ of the group $H_n$.

The other concept is the area of a surface. It is better to call this
area the Toledo invariant of a representation. The remarkable results
of W.~M.~Goldman [Gol1, Corollary C] and D.~Toledo [Tol] say that a
representation is faithful and discrete if (and only if, in the case of
the classic hyperbolic geometry) the `area' of the representation is
`maximal.' In literature (see, for instance, [BIW] and [KMa]), there
are several proofs of Toledo's theorem and neither of them is simple.

First, we study hyperelliptic surfaces. We prove the analog of
W.~M.~Goldman's theorem for hyperelliptic surfaces (Theorem 3.15). The
Teichm\"uller space $\Cal TH_n$ turns out to be supplied with natural
coordinates: The space $\Cal TH_n$ can be described as the space of all
$(2n-6)$-tuples $(z_1,z_2,\dots,z_{2n-6})$ of distinct points on the
ideal boundary $\partial\Bbb D$ that appear in the cyclic order
$z_1,z_2,\dots,z_{2n-6}$ when running once over $\partial\Bbb D$
(Corollary 3.17). These coordinates are natural in the sense that they
have a clear geometric nature and are not related to any arbitrary
choice. Also, they are effective and easily calculable. Besides,
following these ideas, we arrive at a simple and effective criterion
allowing to verify that a given representation is faithful and
discrete. It is worthwhile mentioning a curious fact (we did not find
it in literature) : Every pentagon, i.e., every\footnote{We interpret
as $\PU(1,1)$ the group of all orientation-preserving isometries of
$\Bbb D$.}
representation $\varrho:H_5\to\PU(1,1)$ such that $\varrho(r_i)\ne1$,
is~faithful and discrete (Corollary 3.16). (A complex hyperbolic
version of this fact is discussed in [Ana1, Conjecture 1.2].)

Next, we show that the earthquake group $\Cal E_n$ acts transitively on
$\Cal TH_n$ (Theorem 4.5). This fact can be considered as a continuous
analog\footnote{Maxim Kontsevich convinced us that $\Cal E_n$ is not
finite-dimensional modulo the kernel of its action on $\Cal TH_n$.}
of the well-known Dehn theorem saying that the mapping class group can
be generated by the Dehn twists. (The Dehn twists we use are `integer'
SEs.) Then we prove a discrete variant of Theorem 4.5 --- a sort of the
Dehn theorem: The subgroup of index $2$ in $\Aut H_n$ is generated by
the `integer' SEs (Theorem 4.6).

Finally, we prove W.~M.~Goldman's theorem [Gol1, Corollary C] in
general case (Theorem 5.1). The~idea of the proof is reflected by the
title of this article. We pretend to view a general Riemann surface
$\Sigma$ as if it were a hyperelliptic one and, with a certain
precaution, apply to $\Sigma$ the methods developed in the previous
sections. As in the hyperelliptic case, we establish an effective and
simple criterion of discreteness of a representation of
$G_n:=\pi_1\Sigma$ that involves the construction of a natural
fundamental domain (Remark 5.10). This fundamental domain allows to
visualize the universal family $\Cal F\to\Cal T_n$ of Riemann surfaces,
where $\Cal T_n$ denotes the classic Teichm\"uller space: $G_n$ acting
fibrewise on the trivial bundle $\Bbb D\times\Cal T_n\to\Cal T_n$
provides $\Cal F=\Bbb D\times\Cal T_n/G_n$. The union of the natural
fundamental domains over all fibres is a fundamental domain for the
action of $G_n$ on $\Bbb D\times\Cal T_n$. Yet, we cannot describe
$\Cal T_n$ as explicitely as $\Cal TH_n$. Nevertheless, it is easy to
extend the action of $\Cal E_n$ to $\Cal T_n$ (see Remark 5.24).

Our way of proving the discreteness of a representation, where SEs are
extensively used, resembles a kind of hidden Maskit combination
theorems [Mas]. We think that there is no satisfactory complex
hyperbolic analog of these theorems. The reason is that it is quite
difficult to deduce the discreteness of a `cocompact' group from the
discreteness of its `noncocompact' subgroups appearing after cutting
the corresponding manifold. In our approach, we escape passing to
`noncocompact' groups.

As expected, the complex hyperbolic Toledo theorem [Tol] can be easily
proven (see [Ana2]) by~literally repeating the arguments presented in
this article. Another (unexpected) consequence of our methods is the
fact that $\Cal T_n$ is fibred twice over $\Cal TH_n\subset\Cal T_n$.
Moreover, every point in $\Cal T_n$ is uniquely determined by its
projections to $\Cal TH_n$ [Ana2].

\medskip

{\bf Acknowledgements.} We are very grateful to Fedor Bogomolov, Pedro
Walmsley Frejlich, Carlos Henrique Grossi Ferreira, Nikolay Gusevskii,
and Maxim Kontsevich for their interest to our work.

\bigskip

\centerline{\bf2.~Preliminaries}

\medskip

In our notation, we follow [AGr], except that, for the sake of
convenience, we change the hermitian metric in order to have the
curvature $-1$.

Let $W$ be a two-dimensional $\Bbb C$-vector space equipped with a
hermitian form of signature $+-$. For~a nonisotropic $p\in\Bbb{CP}W$,
define a hermitian form in
$\T_p\Bbb{CP}W\simeq\langle-,p\rangle p^\perp$ as
$\langle t_1,t_2\rangle:=-4\langle p,p\rangle\langle v_1,v_2\rangle$,
where $t_1,t_2\in\T_p{\Bbb C}{\Bbb P}W$, $t_i=\langle-,p\rangle v_i$,
and $v_i\in p^\perp$. The set $\B W$ of negative points in $\Bbb{CP}W$
is simply the open Poincar\'e disc. The set $\overline\B W$ of
nonpositive and the set $\S W$ of isotropic points in $\Bbb{CP}W$ form
the closed Poincar\'e disc and its boundary; all geometrical objects we
deal with live in $\overline\B W$. For distinct
$p_1,p_2\in\overline\B W$, denote by $\G[p_1,p_2]$, $\G(p_1,p_2)$,
$\G(p_1,p_2]$, $\G{\prec}p_1,p_2{\succ}$, etc. the geodesic segments
oriented from $p_1$ to $p_2$ : closed, open, semiopen, full geodesic,
etc.

\smallskip

Let $\Bbb B^2$ denote a closed disc and let
$\varphi:\Bbb B^2\to\overline\B W$ be a piecewise smooth map such that
$\varphi(\partial\Bbb B^2)$ is the union of a finite number of
geodesics and such that $\varphi^{-1}(\S W)\subset\partial\Bbb B^2$ is
finite. Clearly,
$\int\limits_\varphi\omega=\int\limits_{\partial\varphi}P$, where
$\omega$ and $P$ stand for the K\"ahler form and its potential. In
particular, for $p_1,p_2,p_3\in\overline\B W$, the~oriented area of the
triangle $\Delta(p_1,p_2,p_3)$ is given by\footnote{The function $\arg$
takes values in $[-\pi,\pi]$. In the presented formula, the values of
$\arg$ lie in fact in $[-\frac\pi2,\frac\pi2]$.}
$$\Area\Delta(p_1,p_2,p_3)=2\arg\big(-\langle p_1,p_2\rangle\langle
p_2,p_3\rangle\langle p_3,p_1\rangle\big)\leqno{\bold{(2.1)}}$$
(see, for instance, [Gol2] or [AGr, Subsection 5.9]). This formula
works for triangles having no coinciding isotropic vertices. Obviously,
the area of $\Delta(p,p,q)$ vanishes for isotropic $p$. Thus,
$\Area\Delta(p_1,p_2,p_3)$ is continuous while $p_1,p_2,p_3$ run over
$\overline\B W$, assuming different isotropic vertices not to coincide
during the deformation.

Integrating a K\"ahler potential over a closed piecewise geodesic path
$C$ (not necessarily simple), we~obtain the `area' of the `polygon
limited by $C$.' In order to express this area in explicit terms, take
an arbitrary `centre' $c\in\overline\B W$. Let $p_1,p_2,\dots,p_n$ be
successive vertices of $C$. Define
$$\Area(c;C):=\Area(c;p_1,p_2,\dots,p_n):=
\sum_{i=1}^n\Area\Delta(c,p_i,p_{i+1})\leqno{\bold{(2.2)}}$$
(the indices are modulo $n$). Intuitively, this area does not depend on
the choice of $c$. We prefer to give a formal proof of this fact since
it can be useful when we will deal with other invariants different from
the Toledo one.

\medskip

{\bf2.3.~Remark.} For arbitrary
$c,p,q,p_1,p_2,\dots,p_k,q_1,q_2,\dots,q_l\in\overline\B W$, we have
$$\Area(c;p,p_1,p_2,\dots,p_k,q,q_1,q_2,\dots,q_l)=
\Area(c;p,p_1,p_2,\dots,p_k,q)+\Area(c;q,q_1,q_2,\dots,q_l,p)$$
because $\Area\Delta(c,q,p)+\Area\Delta(c,p,q)=0$.

\medskip

So, in order to prove that (2.2) is independent of $c$, we can assume
$n=3$ and the $p_i$'s pairwise distinct. Now, it follows from (2.1)
that
$$\Area\Delta(c;p_1,p_2,p_3)\equiv2\arg\big(\langle
p_1,p_2\rangle\langle p_2,p_3\rangle\langle
p_3,p_1\rangle\big)\mod2\pi$$
for $c$ different from the isotropic $p_i$'s. For such $c$, the
independence follows from the continuity of the triangle area. It is
immediate that $\Area(c;p_1,p_2,p_3)=\Area\Delta(p_1,p_2,p_3)$ for
$c=p_i$. Therefore, it remains to observe that
$\Area(c;p_1,p_2,p_3)=\Area\Delta(p_1,p_2,p_3)$ for $c$ isotropic and
the $p_i$'s pairwise distinct and isotropic, which is straightforward.

\smallskip

For $n\ge5$, let $H_n$ denote the group generated by
$r_1,r_2,\dots,r_n$ with the defining relations $r_i^2=1$,
$i=1,\dots,n$, and $r_n\dots r_2r_1=1$. For even $n$, there is a unique
fully characteristic torsion-free subgroup $G_n$ of index $2$ in $H_n$.
It is constituted by the words of even length in $r_i$'s. As is well
known (see also Proposition 4.1), $G_n$ is the fundamental group of a
closed orientable Riemann surface of genus $\frac n2-1$. For odd $n$,
there is a torsion-free subgroup $T_n$ of index $4$ in $H_n$ which is
the fundamental group of a closed orientable Riemann surface of genus
$n-3$ (see, for instance, [AGG, Subsection 2.1]).

Let $\Cal L:=\PU W$ denote the Lie group of all orientation-preserving
isometries of $\B W$. Denote by $\Cal RH_n$ and $\Cal RG_n$ the spaces
of faithful discrete representations of $H_n$ and of $G_n$ into
$\Cal L$, respectively. The~spaces
$\Cal H_n:=\Cal TH_n:=\Cal RH_n/\Cal L$ and
$\Cal T_n:=\Cal TG_n:=\Cal RG_n/\Cal L$ are the Teichm\"uller spaces of
the groups $H_n$ and $G_n$, i.e., the spaces of conjugacy classes of
the above representations. Each of the two connected components
$\Cal T^-_n$ and $\Cal T^+_n$ of $\Cal T_n$ can be interpreted as the
classic Teichm\"uller space. (The latter appears if we take for
$\Cal L$ the Lie group of all isometries of $\B W$.) Similarly, we
introduce $\Cal H_n^\pm$. The part of $\Cal T^\pm_n$ corresponding to
hyperelliptic surfaces possesses infinitely many connected components
[Mac] which are copies of $\Cal H^\pm_n$ provided by the action of the
mapping class group.

\smallskip

It is easy to see that the involutions in $\Cal L$ are exactly the
reflections in points in $\B W$. Explicitly, in~terms of $\SU W$, such
a reflection $R(q)$ is given by
$R(q):x\mapsto i\Big(x-2\displaystyle\frac{\langle x,q\rangle}{\langle
q,q\rangle}q\Big)$,
$q\in\B W$, $i^2=-1$. Note that $R(q)R(q)=-1$.

\bigskip

\centerline{\bf3.~Hyperelliptic Teichm\"uller Space}

\medskip

Let $\varrho:H_n\to\Cal L$ be a representation. For an arbitrary
$p\in\overline\B W$, define
$$\Area(p;\varrho):=\Area(p_1,p_2,\dots,p_n),\leqno{\bold{(3.1)}}$$
where $p_0:=p$ and $p_i:=\varrho(r_i)p_{i-1}$ (the indices are modulo
$n$). Clearly, we can also define the $p_i$'s starting from
$p=p_j\in\overline\B W$ for an arbitrary $j$ instead of $j=n$.

\medskip

{\bf3.2.~Lemma.} {\sl$\Area(p;\varrho)$ is independent of the choice
of\/ $p$. If\/ $\varrho(r_i)\ne1$ for all\/ $i$, then\/
$\Area(p;\varrho)\equiv n\pi\mod2\pi$.}

\medskip

{\bf Proof.} Without loss of generality, we can assume that
$\varrho(r_i)\ne1$ for all $i$ and choose a representative
$R(q_i)\in\SU W$, $q_i\in\B W$, for every $\varrho(r_i)$. Hence,
$p_i\ne p_{i+1}$ if $p$ is isotropic. It follows from the defining
relations of $H_n$ that $R(q_n)\dots R(q_1)=\varepsilon$, where
$\varepsilon=\pm1$. Take representatives $p_i\in W$ so that
$p_i=R(q_i)p_{i-1}$. In particular, $p_{i+n}=\varepsilon p_i$. It
follows from (2.1) that
$$\Area(p;\varrho)\equiv2\arg\big(\langle p_1,p_2\rangle\langle
p_2,p_3\rangle\dots\langle p_n,p_{n+1}\rangle\big)\mod2\pi.$$
Since $R(q_i)\in\SU W$ and $R(q_i)R(q_i)=-1$, we obtain
$$0\ne\langle p_i,p_{i+1}\rangle=\big\langle
p_i,R(q_{i+1})p_i\big\rangle=\big\langle
R(q_{i+1})p_i,R(q_{i+1})R(q_{i+1})p_i\big\rangle=-\langle
p_{i+1},p_i\rangle.$$
So, $\Area(p;\varrho)\equiv2\arg i^n\equiv n\pi\mod2\pi$, being
$\Area(p;\varrho)$ continuous in $p$
$_\blacksquare$

\medskip

{\bf3.3.~Remark.} For a given representation $\varrho:H_n\to\Cal L$,
define $\varrho J:H_n\to\Cal L$ by
$\varrho J(r_i):=\varrho(r_{n-i})$.\break

\vskip5pt

\noindent
\hskip316pt$\vcenter{\hbox{\epsfbox{Picture.2}}}$

\vskip-117pt

\rightskip145pt

\noindent
Obviously, $\Area\varrho J=-\Area\varrho$. In other words, changing the
cyclic order of the generators alters the sign of the area.

\medskip

In the sequel, we assume without loss of generality that
$\Area\varrho\ge0$.

\medskip

{\bf3.4.~Remark.} Let $p_{i-2}\in\overline\B W$ be a fixed point of
$\varrho(r_ir_{i-1})$. Then, by taking $c=p_{i-2}=p_i$, we can see that

\vskip12pt

\hskip45pt$\Area\Delta(c,p_{i-3},p_{i-2})=
\Area\Delta(c,p_{i-2},p_{i-1})=$

\vskip12pt

\hskip46pt$=\Area\Delta(c,p_{i-1},p_i)=\Area\Delta(c,p_i,p_{i+1})=0$

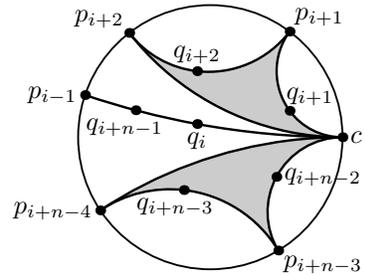
\rightskip0pt

\vskip13pt

\noindent
and, hence, $\Area\varrho\le(n-4)\pi$. When $\Area\varrho=(n-4)\pi$, we
say that $\Area\varrho$ is {\it maximal.}

If $\varrho(r_i)=1$, then $\Area\varrho\le(n-5)\pi$ : `excluding' the
generator $r_i$ we deal in fact with a representation of $H_{n-1}$.

Analogously, if $\varrho(r_ir_{i-1})=1$, then `excluding' the
generators $r_{i-1}$ and $r_i$, we arrive at the representation
$\varrho':H_{n-2}\to\Cal L$. Note that $\Area\varrho=\Area\varrho'$
since $p_{i-2}=p_i$ and
$\Area\Delta(c,p_{i-2},p_{i-1})=-\Area\Delta(c,p_{i-1},p_i)$.
Therefore, $\Area\varrho\le(n-6)\pi$ in this case.

\newpage

\quad

\vskip-6pt

\noindent
\hskip349pt$\vcenter{\hbox{\epsfbox{Picture.3}}}$

\vskip-30pt

\noindent
\hskip345pt$\vcenter{\hbox{\epsfbox{Picture.4}}}$

\rightskip110pt

\vskip-248pt

\medskip

{\bf3.5.~Definition.} Let $q_1,q_2\in\B W$ be distinct. Clearly,
$h^2:=R(q_2)R(q_1)$ for some hyperbolic $h\in\Cal L$. It is easy to see
that $R(h^tq_k)=h^tR(q_k)h^{-t}$, $k=1,2$, and that
$R(q_2)R(q_1)=R(h^tq_2)R(h^tq_1)$ for every $t\in\Bbb R$.

Let $\varrho:H_n\to\Cal L$ be a representation such that
$h^2:=\varrho(r_ir_{i-1})$ is hyperbolic. For every $t\in\Bbb R$,
define a representation $\varrho E_i(t):H_n\to\Cal L$ as follows:
$\varrho E_i(t)(r_j):=\varrho(r_j)$ if $j\notin\{i-1,i\}$ and
$\varrho E_i(t)(r_j):=h^t\varrho(r_j)h^{-t}$, otherwise.

This defines a partial right action of the group $(\Bbb R,+)$ on
representations. We~call $E_i(t)$ a {\it simple earthquake involving\/}
$q_{i-1},q_i$ (SE for short), where $\varrho(r_j)=R(q_j)$,
$j=1,2,\dots n$. Denote by $E_i:=E_i(1)$ the {\it Dehn twist
involving\/} $q_{i-1},q_i$ (DT for short).

\medskip

{\bf3.6~Definition.} If a cycle of isotropic points
$p_1,p_2,\dots,p_k\in\S W$, $k\ge3$, is~listed in the counterclockwise
(clockwise) sense (in particular, the points have to be pairwise
distinct), the cycle is said to be {\it positive\/} ({\it negative\/}).

\medskip

{\bf3.7.~Remark.} Given $p_1,p_2,q_1,q_2\in\S W$, the cycle
$p_1,q_1,p_2,q_2$ is positive or negative if and only if $\G(p_1,p_2)$
and $\G(q_1,q_2)$ intersect in a single point.

If the cycles $p_1,p_2,\dots,p_k\in\S W$, $k\ge3$, and
$p_k,p_{k+1},p_1\in\S W$ are positive, then the cycle
$p_1,p_2,\dots,p_k,p_{k+1}$ is positive.

\medskip

{\bf3.8.~Remark.} Suppose that $\varrho(r_ir_{i-1})$ is hyperbolic.
Then $\Area\varrho=\Area\allowmathbreak\varrho E_i(t)$. Indeed, taking
for $p_{i-2}$ a fixed point of $\varrho(r_ir_{i-1})$, we can see
that\break

\rightskip0pt

\vskip-12pt

\noindent
the $p_j$'s are independent of $t$ and so is $\Area\varrho E_i(t)$.
(See the picture close to Remark 3.4.)

\medskip

{\bf3.9.~Lemma.} {\sl Let\/ $\varrho:H_n\to\Cal L$ be a representation
with maximal\/ $\Area\varrho$. Then, for every\/ $i$, there exists a
suitable\/ $q_i\in\B W$ such that\/ $\varrho(r_i)=R(q_i)$,
$q_{i-1}\ne q_i$, and\/ $\varrho(r_ir_{i-1})$ is hyperbolic. If we take
in\/~{\rm(3.1)}\/ a fixed point of\/ $\varrho(r_ir_{i-1})$ for\/
$p_{i-2}=p_i$, then the cycle\/ $p_i,p_{i+1},\dots,p_{i+n-3}\in\S W$ is
positive.}

\medskip

{\bf Proof.} The first three assertions follow from Remark 3.4 in view
of the fact that the involutions in $\Cal L$ are reflections in points.
As in Remark 3.4, take $c=p_{i-2}=p_i$. The four triangles indicated in
Remark~3.4 are degenerated. Hence, each of the remaining $n-4$ ideal
triangles should have area $+\pi$. In other words, the triangles
$\Delta(c,p_{j-1},p_j)$, $j=i+1,\dots,i+n-3$, are oriented in the
counterclockwise sense. This implies the fourth assertion
$_\blacksquare$

\medskip

{\bf3.10.~Lemma.} {\sl In the situation of Lemma\/ {\rm3.9}, there are
no three collinear points among the\/ $q_j$'s. Moreover,\/
$q_1,q_2,\dots,q_n$ are successive vertices of a convex polygon.}

\medskip

\vskip4pt

\noindent
\hskip319pt$\vcenter{\hbox{\epsfbox{Picture.5}}}$

\rightskip140pt

\vskip-68pt

{\bf Proof.} Suppose that $q_j,q_k,q_l$ are collinear. Acting by $E_j$
or by $E_{j+1}^{-1}$ several times, we can reach a position where
$q_{k-1},q_k,q_l$ are collinear\break

\vskip-1pt

\noindent
$\vcenter{\hbox{\epsfbox{Picture.6}}}$

\leftskip75pt

\vskip-68pt

\noindent
(we diminish $|j-k|>1$). Next, applying $E_l$ or $E_{l+1}^{-1}$ several
times, we arrive at collinear $q_{k-1},q_k,q_{k+1}$. Finally, by means
of some $E_k(t)$, we obtain $q_k=q_{k+1}$. This contradicts Lemma 3.9.

\rightskip0pt

If $q_k$ and $q_l$ are on different sides from
$\G{\prec}q_{j-1},q_j{\succ}$, then $\G{\prec}q_{j-1},q_j{\succ}$ and
$\G[q_k,q_l]$ intersect in some $q\in\B W$. With a suitable $E_j(t)$,
we obtain $q_j=q$, hence, $q_j,q_k,q_l$\break

\vskip-12pt

\leftskip0pt

\noindent
become collinear
$_\blacksquare$

\medskip

{\bf3.11.~Lemma.} {\sl In the situation of Lemma\/ {\rm3.9}, the
points\/ $q_j$, $j\notin\{i-1,i\}$, are on the side of the normal
vector to\/ $\G{\prec}q_{i-1},q_i{\succ}$.}

\newpage

\medskip

\noindent
$\vcenter
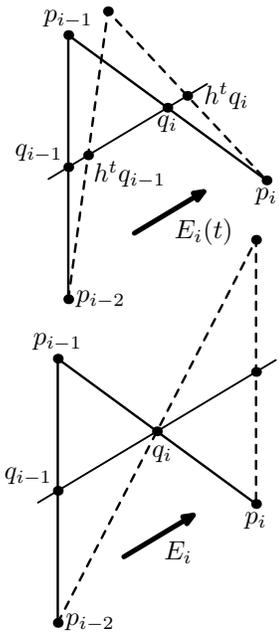
{\hbox{\epsfbox{Picture.7}}}$\hskip16pt$
\vcenter
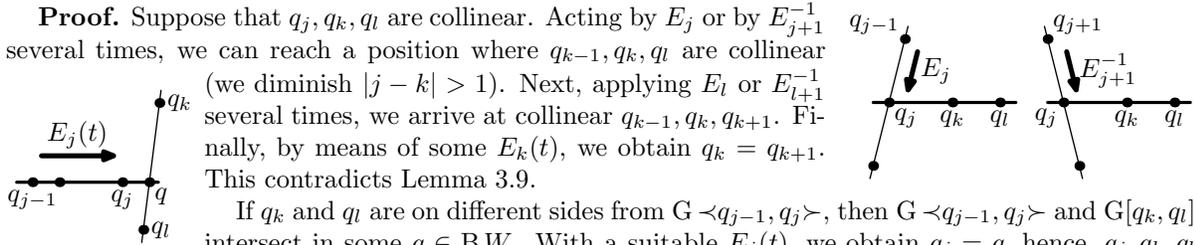
{\hbox{\epsfbox{Picture.8}}}$\hskip16pt$
\vcenter
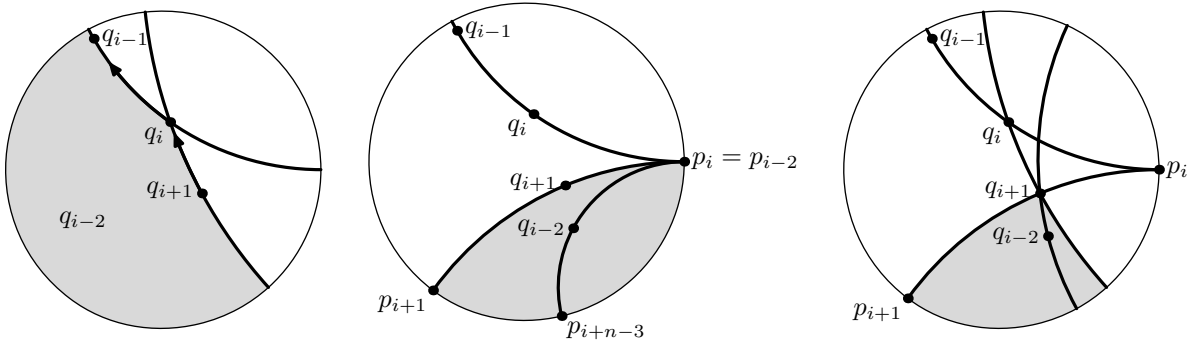
{\hbox{\epsfbox{Picture.9}}}$

\vskip6pt

{\bf Proof.} Due to Lemma 3.10, we can assume all the points $q_j$,
$j\notin\{i-1,i\}$, on the opposite side of the normal vector to
$\G{\prec}q_{i-1},q_i{\succ}$. By Lemma 3.10, this implies that
$q_{i-2}$ is in the region given by the normal vectors to
$\G{\prec}q_i,q_{i-1}{\succ}$ and to $\G{\prec}q_{i+1},q_i{\succ}$,
i.e., in the grey region on the first picture. On the other hand, by
Lemma 3.9, the cycle $p_i,p_{i+1},p_{i+n-3}\in\S W$ is positive, where
$p_{i-2}=p_i\in\S W$ stands for the attractor of $\varrho(r_ir_{i-1})$.
This implies that the geodesic
$\G{\prec}p_{i+n-3},p_{i-2}{\succ}\ni q_{i+n-2}=q_{i-2}$ is entirely on
the side of the normal vector to $\G{\prec}p_i,p_{i+1}{\succ}$ as
illustrated on the second picture. Therefore, the point $q_{i-2}$ is in
the region given by the normal vectors to $\G{\prec}q_{i+1},q_i{\succ}$
and to $\G{\prec}p_i,p_{i+1}{\succ}$ and, thus, the geodesics
$\G{\prec}q_{i-2},q_{i+1}{\succ}$ and $\G{\prec}q_{i-1},q_i{\succ}$
intersect in some point in $\B W$ (see the third picture)
$_\blacksquare$

\medskip

{\bf3.12.~Definition.} Let $\varrho:H_n\to\Cal L$ be a representation
such that $\varrho(r_ir_{i-1})$ is hyperbolic. Denote by $b_i\in\S W$
and by $e_i\in\S W$ the repeller and the attractor of
$\varrho(r_ir_{i-1})$. Put $b_i^i:=b_i$, $e_i^i:=e_i$,
$b_i^j:=\varrho(r_j)b_i^{j-1}$, and $e_i^j:=\varrho(r_j)e_i^{j-1}$. It
follows from the defining relations of $H_n$ that $b_i^{i+n-2}=b_i$ and
$e_i^{i+n-2}=e_i$. We call
$b_i^i,e_i^i,b_i^{i+1},e_i^{i+1},\dots,b_i^{i+n-3},e_i^{i+n-3}\in\S W$
the {\it$i$-cycle\/} of $\varrho$.

\medskip

{\bf3.13.~Proposition.} {\sl Let\/ $\varrho:H_n\to\Cal L$ be a
representation with maximal\/ $\Area\varrho$. Then the\/ $i$-cycle of\/
$\varrho$ is positive.}

\medskip

{\bf Proof.} By Lemma 3.9, the cycles
$b_i^i,b_i^{i+1},\dots,b_i^{i+n-3}$ and
$e_i^i,e_i^{i+1},\dots,e_i^{i+n-3}$ are positive. For suitable points
$q_j\in\B W$, we have $\varrho(r_j)=R(q_j)$. By Lemma 3.11, $q_{i+n-2}$
and $q_{i+1}$ are in the region $D$ given\break

\vskip-7pt

\noindent
\hskip245pt$\vcenter{\hbox{\epsfbox{Picture.10}}}$

\rightskip215pt

\vskip-203pt

\noindent
by the normal vector to $\G{\prec}q_{i-1},q_i{\succ}=\G[b_i^i,e_i^i]$.
So, $e_i^{i+n-3}=R(q_{i+n-2})e_i\in D$. In other words, the cycle
$e_i^{i+n-3},b_i^i,e_i^i$ is positive. Since the geodesics
$\G[e_i^i,e_i^{i+1}]$ and $\G[b_i^i,b_i^{i+1}]$ intersect in
$q_{i+1}\in D\cap\B W$, we have $b_i^{i+1},e_i^{i+1}\in D$ and the
cycle $b_i^i,e_i^i,b_i^{i+1},e_i^{i+1}$ is positive by Remark 3.7. The
fact that the cycles $e_i^i,b_i^{i+1},e_i^{i+1}$ and
$e_i^{i+1},e_i^{i+n-3},e_i^i$ are positive implies that the cycle
$e_i^{i+n-3},e_i^i,b_i^{i+1},e_i^{i+1}$ is positive by Remark 3.7.
Taking into account that the cycle $e_i^{i+n-3},b_i^i,e_i^i$ is
positive, by~Remark 3.7, we get the positive cycle
$e_i^{i+n-3},b_i^i,e_i^i,\allowmathbreak b_i^{i+1},e_i^{i+1}$.

By induction on $j>i$, we can assume that the cycle
$e_i^{i+n-3},b_i^i,e_i^i,\dots,b_i^j,e_i^j$ is positive. The cycle
$e_i^{i+n-3},\allowmathbreak e_i^j,e_i^{j+1}$ is positive. Hence, the
cycle
$e_i^{i+n-3},b_i^i,e_i^i,\dots,\allowmathbreak b_i^j,e_i^j,e_i^{j+1}$
is positive by Remark 3.7. In particular, $b_i^j,e_i^j,e_i^{j+1}$ is
positive. The geodesics $\G[e_i^j,e_i^{j+1}]$ and\break

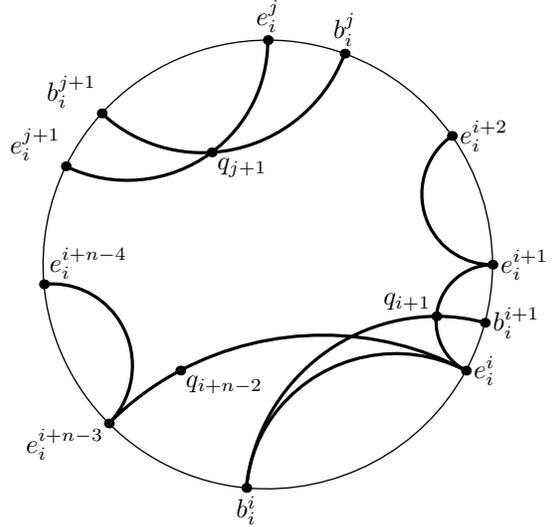
\rightskip0pt

\newpage

\noindent
$\G[b_i^j,b_i^{j+1}]$ intersect (in $q_{j+1}\in\B W$). By Remark 3.7,
the cycle $b_i^j,e_i^j,b_i^{j+1},e_i^{j+1}$ is positive or negative.
Knowing that the cycle $b_i^j,e_i^j,e_i^{j+1}$ is positive, we infer
that $b_i^j,e_i^j,b_i^{j+1},e_i^{j+1}$ is positive and imply that
$e_i^j,b_i^{j+1},e_i^{j+1}$ is positive. Since
$e_i^{i+n-3},b_i^i,e_i^i,\dots,b_i^j,e_i^j,e_i^{j+1}$ is positive,
$e_i^{i+n-3},b_i^i,e_i^i,\dots,b_i^j,e_i^j,b_i^{j+1},e_i^{j+1}$ is
positive by Remark 3.7
$_\blacksquare$

\medskip

{\bf3.14.~Proposition.} {\sl Let\/ $\varrho:H_n\to\Cal L$ be a
representation with hyperbolic\/ $\varrho(r_ir_{i-1})$. If the\/
$i$-cycle of\/ $\varrho$ is positive, then\/ $\varrho\in\Cal RH_n$.}

\medskip

{\bf Proof.} Taking $b_i^i$ for $p_{i-2}$ in (3.1), we obtain the
points $p_{i-2},p_{i-1},\dots,p_{i+n-3}$ which are in fact the points
$b_i^i,e_i^i,b_i^i,b_i^{i+1},b_i^{i+2},\dots,b_i^{i+n-3}$. Since the
$i$-cycle of $\varrho$ is positive, the cycle
$b_i^i,b_i^{i+1},b_i^{i+2},\dots,b_i^{i+n-3}$ is positive and we
conclude that $\Area\varrho=(n-4)\pi$.

Following the natural orientation of $\S W$, we draw an arc
$a_j\subset\S W$ from $b_i^j$ to $e_i^j$ for every
$j=i,i+1,\dots,i+n-3$. The arcs $a_j$ are pairwise disjoint because the
$i$-cycle is positive. We take an arbitrary $p_{i-1}\in\G(q_{i-1},q_i)$
and generate the points $p_j:=\varrho(r_j)p_{j-1}$ so that
$p_{i+n-2},p_{i-1},p_i\in\G_i$, where $G_j:=\G[b_i^j,e_i^j]$. We claim
that $p_{i-1},p_i,\dots p_{i+n-2}$ are the successive vertices of a
convex geodesic\break

\vskip-2pt

\noindent
\hskip39pt$\vcenter
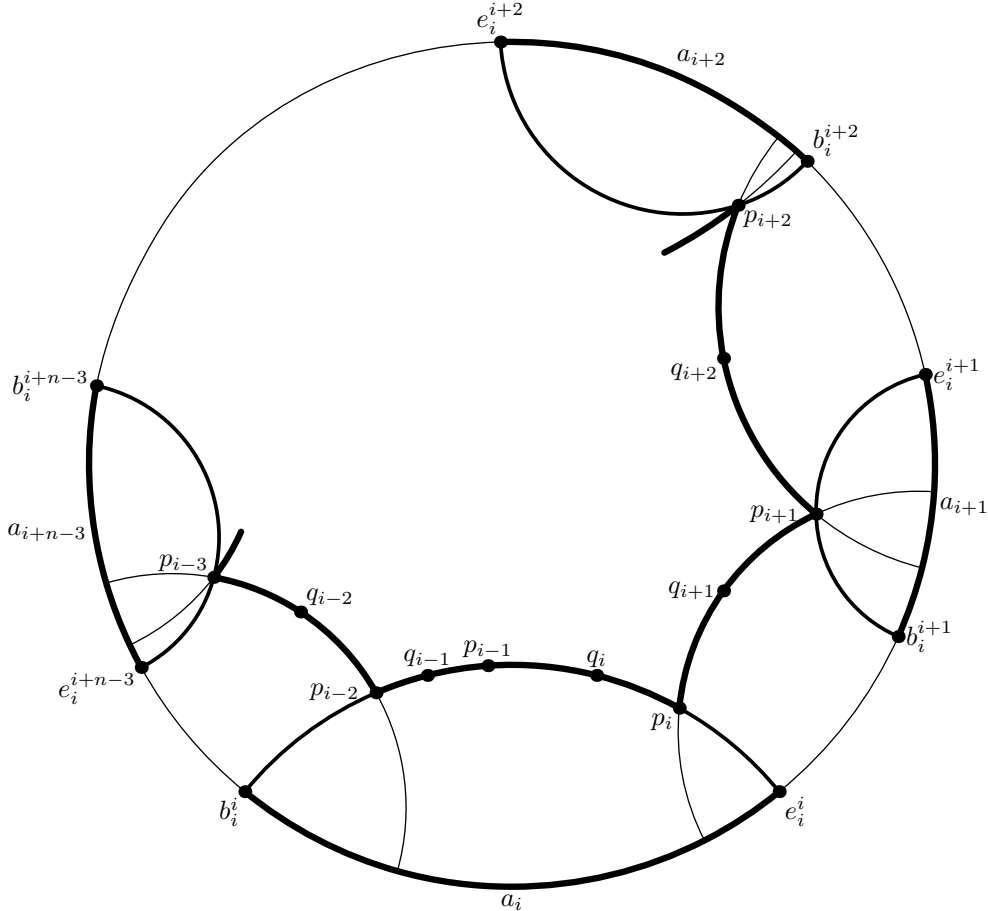
{\hbox{\epsfbox{Picture.11}}}$

\vskip3pt

\noindent
$n$-gon $P_n$. Indeed, $p_j\in\G_j$ for $j=i,i+1,\dots,i+n-3$ because
$G_{j+1}=R(q_{j+1})G_{j}$. For such $j$'s, the vertices of the geodesic
$\Gamma_{j+1}:=\G{\prec}p_j,p_{j+1}{\succ}$ belong to $a_j$ and
$a_{j+1}$ (by convention, $a_{i+n-2}:=a_i$). Hence, $\Gamma_j$ and
$\Gamma_{j+1}$ intersect in $p_j$ and these are the only intersections
between the $\Gamma_j$'s. Since
$\Area\varrho=\Area(p_i,p_{i+1},\dots,p_{i+n-1})=\Area P_n$, the sum of
the interior angles of $P_n$ equals $(n-2)\pi-\Area P_n=2\pi$. By
Poincar\'e's Polyhedron Theorem, $P_n$ is a fundamental polygon for the
group generated by $\varrho(r_j)$ (it~has one cycle of vertices) and
$\varrho$ is faithful and discrete
$_\blacksquare$

\medskip

{\bf3.15.~Theorem.} {\sl Let\/ $\varrho:H_n\to\Cal L$ be a
representation. Then the following statements are equivalent\/{\rm:}

\noindent
$\bullet$~$\varrho\in\Cal RH_n$,\quad
$\bullet$~$\Area\varrho=\pm(n-4)\pi$,\quad $\bullet$~the\/ $i$-cycle
of\/ $\varrho$ is positive or negative.}

\medskip

{\bf Proof} explores standard arguments. We will deal with even $n$
(similar arguments work for odd $n$). Let $\varrho\in\Cal RH_n$.
Clearly, $\varrho|_{_{G_n}}\!\!\in\Cal RG_n$. By definition,
$\Area\varrho=\Area(p_1,p_2,\dots,p_n)$, where
$p_j=\varrho(r_j)p_{j-1}$ for suitable $p_j\in\B W$.

Let $P_n$ be a simple geodesic polygon such that the sum of its
interior angles equals $2\pi$ and let $v_1,v_2,\dots,v_n$ stand for the
successive vertices of $P_n$ listed in the counterclockwise sense. Let
$q_j$ denote the middle point of $\G[v_{j-1},v_j]$. By Poincar\'e's
Polyhedron Theorem, $P_n$ is a fundamental polygon for the group
generated by $R(q_j)$ and, thus, we arrive at some
$\varrho_0\in\Cal RH_n$.

Let us define a continuous $H_n$-equivariant map $\varphi:\B W\to\B W$.
Put $\varphi v_j=p_j$ and define $\varphi$ linearly on the geodesic
$\G[v_{j-1},v_j]$; so, $\varphi\G[v_{j-1},v_j]=\G[p_{j-1},p_j]$. Next,
extend $\varphi$ continuously to $\varphi:P_n\to\B W$. Finally, put
$\varphi\big(\varrho_0(h)p\big)=\varrho(h)\varphi(p)$ for all
$h\in H_n$ and $p\in P_n$. The map $\varphi$ induces a continuous
map\break

\vskip-3pt

\noindent
\hskip380pt$\vcenter{\hbox{\epsfbox{Picture.12}}}$

\rightskip80pt

\vskip-115pt

\noindent
$\psi:\Sigma_0\to\Sigma$, where $\Sigma_0:=\B W/\varrho_0G_n$ and
$\Sigma:=\B W/\varrho G_n$ are Riemann surfaces of genus $\frac n2-1$.
By construction, $\pi_1\psi:\pi_1\Sigma_0\to\pi_1\Sigma$ is an
isomorphism, hence,
$H_2\psi:H_2(\Sigma_0,\Bbb Z)\to H_2(\Sigma,\Bbb Z)$ is an isomorphism
and
$\int\limits_\psi\omega'=\mp\Area\Sigma=\pm2\pi\chi(\Sigma)=
\mp2(n-4)\pi$,
where $\omega'$ stands for the K\"ahler form of $\Sigma$. On the other
hand, $P_n\cup\varrho_0(r_i)P_n$ is a fundamental polygon for
$\varrho_0G_n$, therefore,
$\int\limits_{\psi}\omega'=2\int\limits_{\varphi|_{_{P_n}}}\omega=
2\int\limits_{\varphi|_{_{\partial
P_n}}}P=-2\Area(p_1,p_2,\dots,p_n)=-2\Area\varrho$,
where $P$ stands for a K\"ahler potential of $\B W$. Consequently,
$\Area\varrho=\pm(n-4)\pi$
$_\blacksquare$

\medskip

\noindent
$\vcenter{\hbox{\epsfbox{Picture.13}}}$

\leftskip80pt

\vskip-66pt

From Remark 3.3 and Lemma 3.2, we obtain the

\rightskip0pt

\medskip

{\bf3.16.~Corollary.} {\sl Let\/ $\varrho:H_5\to\Cal L$ be a
representation such that\/ $\varrho(r_i)\ne1$ for all\/~$i$. Then\/
$\varrho\in\Cal RH_5$
$_\blacksquare$}

\medskip

Note that Theorem 3.15 provides an effective criterion of discreteness:
In order to verify that some $q_1,\dots,q_n\in\B W$ subject to the
relation $R(q_n)\dots R(q_1)=\pm1$\break

\leftskip0pt

\vskip-12pt

\noindent
provide a representation $\varrho\in\Cal RH_n$, we can explicitly find
the $b_i^j$'s and $e_i^j$'s and check if the $i$-cycle of $\varrho$ is
positive or negative.

Also, Theorem 3.15 yields some explicit description of the two
components $\Cal H_n^+$ and $\Cal H_n^-$ (related to the sign of
$\Area\varrho$) of $\Cal H_n$. Let
$\Bbb S^1_+:=\big\{z\in\Bbb C\mid |z|=1,\Im z>0\big\}$ and let
$$\Cal K_n^+:=\big\{(z_1,z_2,\dots,z_{2n-7},z_{2n-6})\mid z_j\in\Bbb
S^1_+,\text{ \rm the cycle }z_1,z_2,\dots,z_{2n-7},z_{2n-6}\text{ \rm
is positive}\big\}$$
(similarly, we define $\Bbb S^1_-$ and $\Cal K_n^-$).

\noindent
\hskip74pt$\vcenter
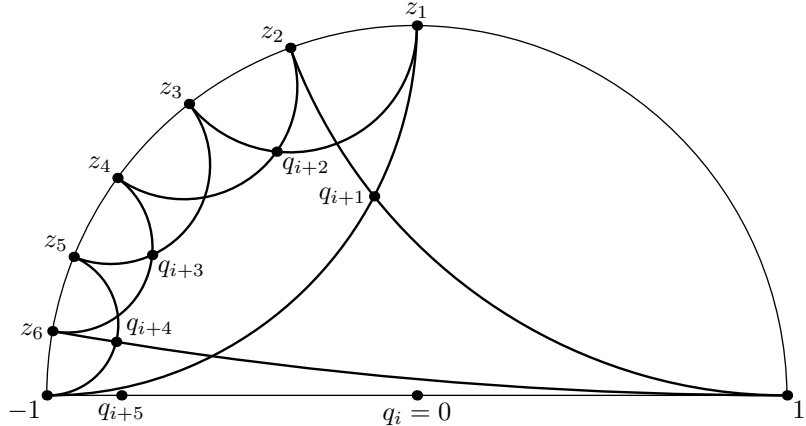
{\hbox{\epsfbox{Picture.14}}}$

\medskip

{\bf3.17.~Corollary.} {\sl$\Cal H_n^\pm\simeq\Cal K_n^\pm$.}

\medskip

{\bf Proof.} Identify $\overline\B W$ with the unitary disc
$\big\{z\in\Bbb C\mid|z|\le1\big\}$. Let $[\varrho]\in\Cal H_n^+$.
Conjugating $\varrho$ with an element in $\Cal L$, we can assume that
$b_i^i=-1$, $e_i^i=1$, and $q_i=0$. This provides
$$(z_1,z_2,\dots,z_{2n-7},z_{2n-6}):=
(b_i^{i+1},e_i^{i+1},\dots,b_i^{i+n-3},e_i^{i+n-3}).$$
In other words, we obtain a map $\Cal H_n^+\to\Cal K_n^+$.

Conversely, for given $(z_1,z_2,\dots,z_{2n-7},z_{2n-6})\in\Cal K_n^+$,
define $q_i:=0$, $q_{i+1}:=\G[-1,z_1]\cap\G[1,z_2]$,
$q_{i+k}:=\G[z_{2k-3},z_{2k-1}]\cap\G[z_{2k-2},z_{2k}]$ for
$k=2,3,\dots,n-3$, and $q_{i+n-2}:=\G[z_{2n-7},-1]\cap\G[z_{2n-6},1]$.
It is easy to see that the isometry
$h:=R(q_{i+n-2})\dots R(q_{i+2})R(q_{i+1})\in\SU W$ fixes the points
$-1$ and~$1$. If $h=\pm1$, we obtain a representation
$\varrho_0:H_{n-2}\to\Cal L$. Taking $p_i=1$, we arrive at
$\Area\varrho_0=\Area(1,z_2,z_4,\dots,z_{2n-6})=(n-4)\pi$, which
contradicts Remark 3.4. Therefore, $h$ is hyperbolic with the axis
$\G[-1,1]$ and there exists a unique $q_{i+n-1}\in\G(-1,1)$ such that
$h=R(q_{i+n-1})R(q_i)$. In other words,
$R(q_{i+n-1})\dots R(q_{i+1})R(q_i)=\pm1$, providing a representation
$\varrho$ whose $i$-cycle is positive
$_\blacksquare$

\medskip

Note that the indicated identification is effectively calculable with a
simple algorithm. It is easy to show that the points $q_j$ can be
algebraically expressed in terms of the $z_k$'s (not involving
radicals, when using the Klein model). This is why we can treat
$\Cal H_n^{\pm}$ as a `rational variety.'

We are going to study the space $\Cal H_n^{\pm}$ in detail in
subsequent articles. In particular, we would like to describe the
standard hermitian and complex structures of $\Cal H_n^{\pm}$ in terms
of the $z_k$'s. We can introduce a complex structure on $\Cal H_n^+$ by
taking the $q_i$'s, $i=2,\dots,n-2$, as complex coordinates that vary
in the open upper half-disc. Taking $q_1=0$, we can reconstruct
$q_n\in(-1,0)$ and $q_{n-1}$ from given $q_i$'s, $i=2,\dots,n-2$.
However, it is easy to see that the DT $E_3$ is not holomorphic with
respect to this structure. So, it is not the genuine one. (The DT $E_3$
belongs to the hyperelliptic mapping class group (see Section 4) which
is known to be the group of holomorphic automorphisms of $\Cal H_n^+$.)

\bigskip

\centerline{\bf4.~Earthquake Group and Hyperlliptic Mapping Class
Group}

\medskip

Let $\tau:\Sigma\to\Sigma/\iota\simeq\Bbb{CP}^1$ be a hyperelliptic
Riemann surface of genus $g$, where $\iota:\Sigma\to\Sigma$ stands for
the hyperelliptic involution of $\Sigma$. Put $n:=2g+2$ and denote by
$f_1,f_2,\dots,f_n\in\Sigma$ the fixed points of $\iota$. Let
$F\le\Cal L$ stand for the fundamental group of $\Sigma=\B W/F$ and
$\pi:\B W\to\Sigma$, for the universal covering of $\Sigma$. 

\medskip

{\bf4.1.~Proposition {\rm[Mac]}.} {\sl$\Sigma\simeq\B W/\varrho G_n$
for some\/ $\varrho\in\Cal RH_n$. If\/ $\varrho\in\Cal RH_n$, then\/
$\Sigma\simeq\B W/\varrho G_n$ is hyperelliptic.}

\medskip

\noindent
$\CD\B W@>R>>\B W\\@V\pi VV@V\pi VV\\\Sigma@>\iota>>\Sigma\endCD$

\vskip-53pt

\noindent
\hskip323pt$\vcenter{\hbox{\epsfbox{Picture.15}}}$

\leftskip90pt\rightskip136pt

\vskip-76pt

{\bf Proof} explores many well-known arguments.\break
For every $q\in Q:=\pi^{-1}\{f_1,\dots,f_n\}$, there exists a unique
$R\in\Cal L$ inducing in $\Sigma$ the isometry $\iota$ such that
$Rq=q$. Clearly, $F^R=F$. It is easy to see that $R=R(q)$. Indeed, the
isometry $R$ is elliptic\break

\vskip-12pt

\leftskip0pt

\noindent
and $R^2$ induces in $\Sigma$ the isometry $\iota^2=1$. Therefore,
$R^2\in F$, which implies $R^2=1$ because the isometries in  $F$ have
no fixed points in $\B W$.\break

\vskip-12pt

\rightskip0pt

\noindent
For $q_1,q_2\in Q$, the product $R(q_1)R(q_2)$ induces in $\Sigma$ the
isometry $\iota^2=1$. This implies $R(q_1)R(q_2)\in F$.

Choose and fix a point $p\in\B W$ that belongs to no geodesic joining
points in $Q$. Let $q_i$ denote a point in $\pi^{-1}(f_i)$ closest to
$p$, i.e., $\dist(p,q_i)\le\dist(p,fq_i)$ for all $f\in F$. Note that
$q_i$ is also a point in $\pi^{-1}(f_i)$\break

\newpage

\noindent
closest to $R(q_i)p$, i.e.,
$\dist(R(q_i)p,q_i)\le\dist(R(q_i)p,fq_i)$ for all $f\in F$. This
follows from $f^{R(q_i)}\in F$.\break

\vskip-4pt

\noindent
$\vcenter
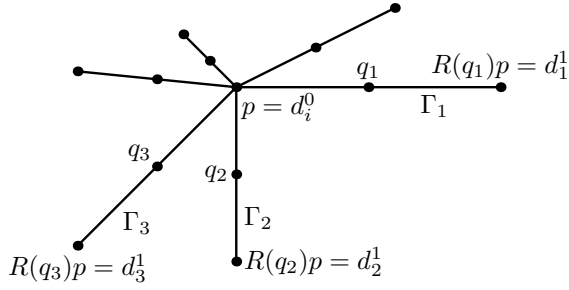
{\hbox{\epsfbox{Picture.16}}}$

\vskip-112pt

\leftskip222pt

Define $d_i^0:=p$, $d_i^1:=R(q_i)p$, $\Gamma_i:=\G[d_i^0,d_i^1]$, and
$\Gamma:=\bigcup\Gamma_i$. Note that
$\dist(p,q_j)\le\dist(d_i^\varepsilon,fq_j)$ for all $f\in F$. Indeed,
this is clear for $\varepsilon=0$. Taking $\varepsilon=1$, we have

\vskip11pt

\noindent
\hskip22pt$\dist(p,q_j)\le\dist\big(p,R(q_i)R(q_j)f^{R(q_j)}q_j\big)=$

\vskip11pt

\noindent
\hskip38pt$=\dist\big(R(q_i)p,fq_j\big)=\dist(d_i^1,fq_j)$

\vskip11pt

\noindent
for all $f\in F$.

\leftskip0pt

\noindent
\hskip397pt$\vcenter{\hbox{\epsfbox{Picture.17}}}$

\rightskip65pt

\vskip-105pt

Clearly, $\pi$ identifies the $d_i^1$'s. Let us show that $\pi$
identifies in $\Gamma$ only the points $d_i^1$. The point $\pi p$ is
not a fixed point of $\iota$, hence, $\pi$ cannot identify $d_i^1$ and
$d_j^0=p$. Suppose that $\pi c_i=\pi c_j$ for $c_i\in\Gamma_i$ and
$c_j\in\Gamma_j$. This means that $fc_j=c_i$ for some $f\in F$. We
assume that $q_i\ne fq_j$ since $i=j$ and $f=1$, otherwise. Let
$d_j^\delta$ denote an end of $\Gamma_j$ closest to $c_j$ and
$d_i^\varepsilon$, an end of $\Gamma_i$ closest to $c_i$. Since
$p\notin\G{\prec}q_i,fq_j{\succ}$ by the choice of $p$, we have
$fq_j\notin\G{\prec}d_i^\varepsilon,q_i{\succ}\ni p$. Without loss of
generality, we can assume that
$\dist(d_i^\varepsilon,c_i)\le\dist(d_j^\delta,c_j)=
\dist(fd_j^\delta,fc_j)$
and that $c_j\ne d_j^\delta$. The lenght\footnote{When a path $x$ ends
with the start point of a path $y$, we denote by $x\cup y$ their
path-product.}
of the path $\gamma:=\G[d_i^\varepsilon,c_i]\cup\G[fc_j,fq_j]$ is less
or equal than that of
$\G[fd_j^\delta,fc_j]\cup\G[fc_j,fq_j]=\G[fd_j^\delta,fq_j]$.
Consequently, $\dist(d_i^\varepsilon,fq_j)\le\dist(fd_j^\delta,fq_j)$.
Since

\rightskip0pt

$$\dist(fd_j^\delta,fq_j)=\dist(d_j^\delta,q_j)=
\dist(p,q_j)\le\dist(d_i^\varepsilon,fq_j),$$
we conclude that $\dist(d_i^\varepsilon,fq_j)=\dist(fd_j^\delta,fq_j)$.
Thus, $\gamma=\G[d_i^\varepsilon,fq_j]$ and, in view of
$fq_j\notin\G{\prec}d_i^\varepsilon,q_i{\succ}$, the point $c_i$ has to
coincide with $d_i^\varepsilon$. From
$\dist(fd_j^\delta,fq_j)=\dist(d_i^\varepsilon,fq_j)$, we conclude that
$fc_j=c_i\allowmathbreak=fd_j^\delta$. A contradiction.

\smallskip

The involution $\iota$ identifies one half of $\pi\Gamma_i$ with the
other since $R(q_i)$ induces in $\Sigma$ the isometry $\iota$. Those
are the only identifications in $\pi\Gamma$ by $\iota$. The curve
$\mu_i:=\tau\pi\Gamma_i$ begins with $\tau\pi p$ and ends with
$\tau f_i$. The only pairwise intersection between the $\mu_i$'s is
$\tau\pi p$. We can assume that $\mu_1,\mu_2,\dots,\mu_n$ are listed in
the clockwise sense with respect to the standard orientation of
$\Bbb{CP}^1$. For every $i$, choose a small open disc
$D_i\subset\Bbb{CP}^1$ centred at $\tau f_i$ such that $D_i$ intersects
$\mu_i$ in some final segment $s_i\subset\mu_i$ and such that the
$D_i$'s are pairwise disjoint. Also, choose a simple closed curve
$\omega_i\subset D_i$ that begins with $p_i\in s_i$, $p_i\ne\tau f_i$,
and winds once around $\tau f_i$ in the clockwise sense. Let
$\sigma_i\subset\mu_i$ denote the segment that begins with $\tau\pi p$
and ends\break

\vskip-5pt

\noindent
\hskip35pt$\vcenter
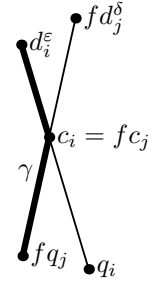
{\hbox{\epsfbox{Picture.18}}}$\hskip35pt$\vcenter
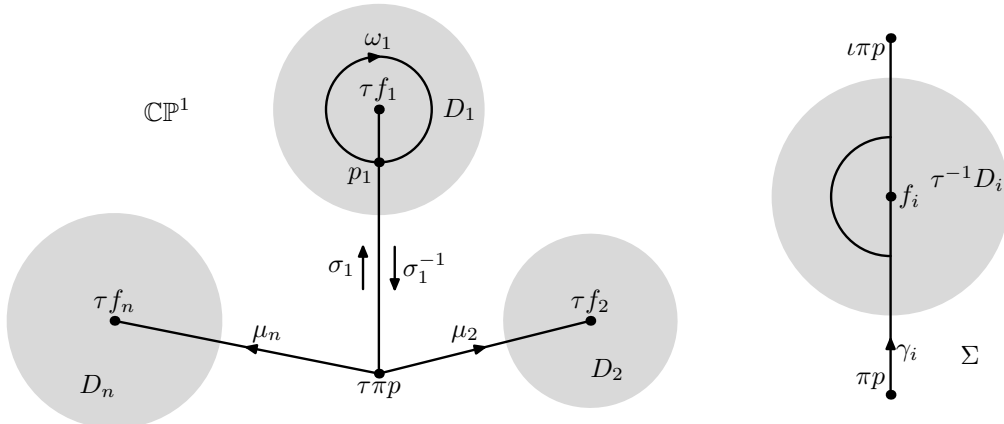
{\hbox{\epsfbox{Picture.19}}}$

\noindent
with $p_i$. The lift $\eta$ of
$\bigcup(\sigma_i\cup\omega_i\cup\sigma^{-1}_i)$ based at $\pi p$ is
contractible in $\Sigma$ and runs over almost all $\pi\Gamma$.
Deforming inside the open discs $\tau^{-1}D_i$ the parts of $\eta$ that
are lifts of the $\omega_i$'s, we arrive at a curve
$\gamma\subset\pi\Gamma$ contractible in $\Sigma$. Clearly, $\gamma$
runs over all $\pi\Gamma$, once over each $\pi\Gamma_i$. For the same
reason, every element in $F=\pi_1(\Sigma,\pi p)$ is represented by a
curve included in $\pi\Gamma$ since the group
$\pi_1\big(\Bbb{CP}^1\setminus\{\tau f_1,\tau f_2,\dots,\tau
f_n\},\tau\pi p\big)$
is generated by the elements $[\sigma_i\cup\omega_i\cup\sigma^{-1}_i]$.
We assume that $\gamma_i:=\pi\Gamma_i$ begins with $\pi p$, passes
through $f_i$, and ends with $\iota\pi p$. The group
$F=\pi_1(\Sigma,\pi p)$ is generated by the elements
$[\gamma_{i-1}\cup\gamma_i^{-1}]$ (the indices are modulo $n$).
Therefore, the elements
$\lambda_{2i}:=\gamma_1\cup\gamma_2^{-1}\cup\dots\cup\gamma_{2i-1}\cup
\gamma_{2i}^{-1}$
and
$\lambda_{2i+1}:=\gamma_2\cup\gamma_3^{-1}\cup\dots\cup\gamma_{2i}\cup
\gamma_{2i+1}^{-1}$,
$i=1,2,\dots,\frac n2$, also generate $F$. Note that $\lambda_n=\gamma$
is contractible in $\Sigma$ by construction.

The lift $\Lambda_{2i}$ of $\lambda_{2i}$ based at $p$ is formed by
$\Gamma_1,\Gamma'_2,\dots,\Gamma'_{2i}$, conjugates of $\Gamma_i$'s.
Let $q'_i\in Q$ denote the middle point of $\Gamma'_i$. Then
$R(q'_{2i})\dots R(q'_2)R(q_1)p$ is the end of $\Lambda_{2i}$. Hence,
$[\lambda_{2i}]=R(q'_{2i})\dots R(q'_2)R(q_1)$. The lift
$\Lambda_{2i+1}$ of $\lambda_{2i+1}$ based at $p$ is formed by
$\Gamma_2,\Gamma''_3,\dots,\Gamma''_{2i+1}$, conjugates of
$\Gamma_i$'s. Let $q''_i\in Q$ denote the\break

\vskip-2pt

\noindent
\hskip265pt$\vcenter{\hbox{\epsfbox{Picture.20}}}$

\rightskip193pt

\vskip-107pt

\noindent
middle point of $\Gamma''_i$. Then $R(q''_{2i+1})\dots R(q''_3)R(q_2)p$
is the end of $\Lambda_{2i+1}$. Hence,
$[\lambda_{2i+1}]=R(q''_{2i+1})\dots R(q''_3)R(q_2)$. We put
$q''_2:=q_2$ and $q'_1:=q_1$.

Let us show that $q''_j=R(q_1)q'_j$. Note that if some $f\in F$ maps a
point in $\Gamma'_i$ to a point in $\Gamma''_i$, then
$\Gamma''_i=f\Gamma'_i$. Since $R(q_1)$ maps the beginning of
$\Gamma'_2$ to $p$, we conclude that $R(q_1)R(q'_2)$ maps the end of
$\Gamma'_2$ to $p$, the beginning of $\Gamma_2$. Hence,
$\Gamma''_2=R(q_1)R(q'_2)\Gamma'_2=R(q_1)\Gamma'_2$. By induction on
$j$, we assume that
$\Gamma''_j=R(q_1)R(q'_j)\Gamma'_j=R(q_1)\Gamma'_j$. Since
$R(q_1)$\break

\vskip-11pt

\rightskip0pt

\noindent
maps the end of $\Gamma'_j$ (which is the beginning of $\Gamma'_{j+1}$)
to the end of $\Gamma''_j$ (which is the beginning of
$\Gamma''_{j+1}$), we conclude that $R(q_1)R(q'_{j+1})$ maps the end of
$\Gamma'_{j+1}$ to the beginning of $\Gamma''_{j+1}$. Hence,
$\Gamma''_{j+1}=R(q_1)R(q'_{j+1})\Gamma'_{j+1}=R(q_1)\Gamma'_{j+1}$.

Consequently, $R(q'_n)\dots R(q'_2)R(q'_1)=1$, which generates a
representation $\varrho:H_n\to\Cal L$ such that $\varrho G_n=F$. Being
$G_n$ and $F$ the fundamental groups of Riemann surfaces of the same
genus, $\varrho|_{_{G_n}}$ is an isomorphism. So is $\varrho$.

The converse can be readily shown with the help of the fundamental
polygon for $\varrho H_n$ constructed in the proof of Proposition 3.14
$_\blacksquare$

\medskip

The formal multiplicative group generated by $n$ copies
$\big\{E_i(t)\mid t\in\Bbb R\big\}$, $i=1,2,\dots,n$, of $(\Bbb R,+)$
is denoted by $\Cal E_n$ and called {\it earthquake group.} We
distinguish the parts $\Cal R^+H_n$ and $\Cal R^-H_n$ of $\Cal RH_n$
related to the sign of the area of a representation. Due to Remark 3.8,
$\Cal E_n$ acts from the right by means of SEs on\break

\vskip-2pt

\noindent
$\vcenter{\hbox{\epsfbox{Picture.22}}}$

\leftskip147pt

\vskip-150pt

\noindent
$\Cal R^\pm H_n$ and, hence, on $\Cal H_n^\pm$. Later (see Remark 5.24)
we will extend this action to $\Cal RG_n$ and to $\Cal T_n$.

\noindent
\hskip258pt$\vcenter
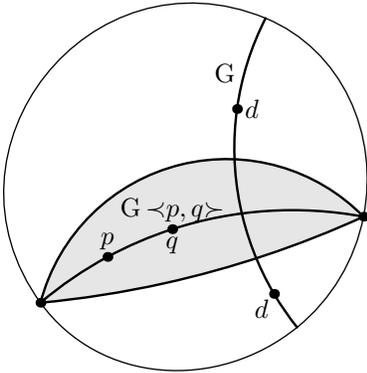
{\hbox{\epsfbox{Picture.21}}}$

\rightskip55pt

\vskip-45pt

{\bf4.2.~Lemma {\rm[Ana1]}.} {\sl Let\/ $\varrho,\varrho'\in\Cal RH_5$
be such that\/ $\varrho(r_4)=\varrho'(r_4)$,
$\varrho(r_5)=\varrho'(r_5)$, and\/ $\Area\varrho=\Area\varrho'$. Then
we can obtain\/ $\varrho'$ from\/ $\varrho$ by means of a finite number
of SEs of the types\/ $E_2(t)$ and\/ $E_3(t)$
$_\blacksquare$}

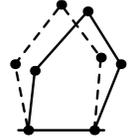
\rightskip0pt

\medskip

{\bf4.3.~Remark.} Let $p,q\in\B W$ be distinct and let $\G$ be a full
geodesic different from $\G{\prec}p,q{\succ}$ and intersecting
$\G{\prec}p,q{\succ}$ in some point in $\B W$. Then, on any side from
$\G{\prec}p,q{\succ}$, there exists some $d\in\G\cap\B W$ such that
$R(d)R(q)R(p)$ is hyperbolic. Indeed, the points $d\in\B W$ making
$R(d)R(q)R(p)$ parabolic form two curves (hypercycles) equidistant from
$\G{\prec}p,q{\succ}$. The isometry $R(d)R(q)R(p)$ is hyperbolic\break

\vskip-12pt

\leftskip0pt

\noindent
exactly when $d$ is outside the band limited by these curves.

\medskip

{\bf4.4.~Lemma.} {\sl$\Cal E_5$ acts transitively on\/
$\Cal R^\pm H_5$.}

\medskip

\newpage

\noindent
\hskip273pt$\vcenter{\hbox{\epsfbox{Picture.23}}}$\hskip7pt$\vcenter
{\hbox{\epsfbox{Picture.24}}}$

\rightskip187pt

\vskip-93pt

{\bf Proof.} Let $\varrho,\varrho'\in\Cal R^+H_5$, that is,
$\Area\varrho=\Area\varrho'=\pi$. For suitable $q_i,q'_i\in\B W$, we
have $\varrho(r_i)=R(q_i)$ and $\varrho'(r_i)=R(q'_i)$.

By Lemma 3.11, the points $q_1,q_2,q_3$ are on the side of the normal
vector to $\G{\prec}q_4,q_5{\succ}$. Let $\G$ be a geodesic passing
through $q'_1$ and intersecting $\G{\prec}q_4,q_5{\succ}$ in some point
in $\B W$. By Remark 4.3, $R(q''_1)R(q_5)R(q_4)$ is hyperbolic for some
$q''_1\in\G\cap\B W$ on the mentioned side. Hence,
$R(q''_1)R(q_5)R(q_4)=$

\rightskip0pt

\noindent
$R(q''_2)R(q''_3)$ for some $q''_2,q''_3\in\B W$, which provides some
$\varrho''\in\Cal RH_5$ by Corollary 3.16. Since $q''_1$ is on the side
of the normal vector to $\G{\prec}q_4,q_5{\succ}$, it follows that
$\Area\varrho''=\pi$ by Lemma 3.11. By Lemma 4.2, after applying to
$\varrho$ a finite number of SEs, we can assume that $q_i=q''_i$,
$i=1,2,3$. Some SE involving $q_4,q_5$ puts $q_5$ into
$\G{\prec}q_1,q'_1{\succ}=\G$. Now, some SE involving $q_5,q_1$
provides $q_1=q'_1$.

\vskip3pt

\noindent
$\vcenter{\hbox{\epsfbox{Picture.25}}}$

\vskip-73pt

\noindent
\hskip362pt$\vcenter{\hbox{\epsfbox{Picture.26}}}$

\rightskip98pt

\leftskip78pt

\vskip-96pt

Since $q_1=q'_1\ne q'_2$, after applying (if necessary) some SE that
involves $q_2,q_3$, we obtain $q'_2\notin\G{\prec}q_1,q_2{\succ}$. By
Remark 4.3, there exist $q''_3\in\G{\prec}q_2,q'_2{\succ}$ and
$q''_4,q''_5\in\B W$ such that the relation
$R(q''_3)R(q_2)R(q_1)=R(q''_4)R(q''_5)$ provides some
$\varrho''\in\Cal RH_5$ with $\Area\varrho''=\pi$. As above, by Lemma
4.2, we can assume that $q_i=q''_i$, $i=3,4,5$. By means of some SE
involving $q_2,q_3$,\break

\vskip-12pt

\leftskip0pt\rightskip0pt

\noindent
we arrive at $q_2=q'_2$. It remains to apply Lemma 4.2 once more
$_\blacksquare$

\medskip

{\bf4.5.~Theorem.} {\sl$\Cal E_n$ acts transitively on\/
$\Cal R^\pm H_n$.}

\medskip

{\bf Proof.} Let $\varrho,\varrho'\in\Cal R^+H_n$, i.e.,
$\Area\varrho=\Area\varrho'=(n-4)\pi$. For suitable $q_i,q'_i\in\B W$,
we have $\varrho(r_i)=R(q_i)$ and $\varrho'(r_i)=R(q'_i)$.

\vskip6pt

\noindent
$\vcenter
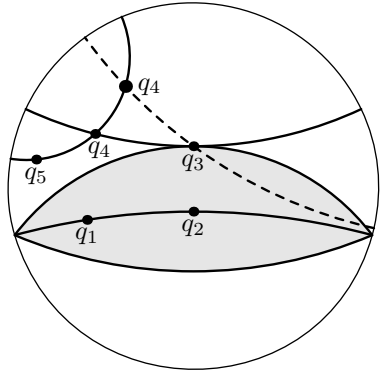
{\hbox{\epsfbox{Picture.27}}}$

\leftskip151pt

\vskip-146pt

The isometry $R(q_3)R(q_2)R(q_1)$ is hyperbolic because
$\varrho\in\Cal RH_n$. Indeed, if it is parabolic, $q_3$ belongs to the
hypercycle
$H=\big\{q\in\B W\mid R(q)R(q_2)R(q_1)\text{ is parabolic}\big\}$.
Applying a `small' SE involving $q_4,q_5$ if necessary, we can assume
$\G{\prec}q_3,q_4{\succ}$ to be transversal to $H$ at $q_3$. Now a
suitable SE involving $q_3,q_4$ provides an elliptic
$R(q_3)R(q_2)R(q_1)$ (see Remark 4.3). A contradiction.

\vskip3pt

\noindent
\hskip253pt$\vcenter{\hbox{\epsfbox{Picture.28}}}$

\rightskip56pt

\vskip-68pt

Hence, there exist $b,d\in\B W$ such that
$R(d)R(b)\allowmathbreak R(q_3)R(q_2)R(q_1)=1$, which generates a
representation $\varrho_0:H_5\to\Cal L$. The relation
$R(q_n)\dots R(q_5)R(q_4)R(b)\allowmathbreak R(d)=1$ generates a
representation $\varrho_1:H_{n-1}\to\Cal L$. Take for $p_0\in\S W$ a
fixed point of $R(b)R(d)=R(q_3)\allowmathbreak R(q_2)R(q_1)$. By Remark
2.3,

\rightskip0pt\leftskip0pt

\vskip-18pt

$$\Area\varrho=\Area(p_0,p_1,p_2,p_3,p_4,\dots,p_{n-1})=
\Area(p_0,p_1,p_2,p_3)+\Area(p_3,p_4,\dots,p_{n-1},p_0)=$$
$$=\Area(p_0,p_1,p_2)+\Area(p_3,p_4,\dots,p_{n-1})=
\Area\varrho_0+\Area\varrho_1$$
due to $p_0=p_3$. By Remark 3.4, $\Area\varrho_0=\pi$ and
$\Area\varrho_1=(n-5)\pi$. By Theorem 3.15, $\varrho_0\in\Cal R^+H_5$
and $\varrho_1\in\Cal R^+H_{n-1}$.

We are going to express every SE of $\varrho_1$ in terms of suitable
SEs of $\varrho$ and an SE involving $d,b$ (the~latter is simply a
rechoice of $b$ and $d$). The SEs of $\varrho_1$ involving
$q_{i-1},q_i$, $i=5,6,\dots,n$, are in fact some SEs of $\varrho$. All
we need is to execute the SEs of $\varrho_1$ involving the pairs
$b,q_4$ and $q_n,d$. By symmetry, we deal only with the first one.

By Remark 4.3, we can find $q''_3\in\G{\prec}b,q_4{\succ}\cap\B W$ on
the side of the normal vector to $\G{\prec}b,d{\succ}$ such that
$R(d)R(b)R(q''_3)$ is hyperbolic and, hence,
$R(d)R(b)R(q''_3)=R(q''_1)R(q''_2)$ for some $q''_1,q''_2\in\B W$.
As~in the proof of Lemma 4.4 (that is, by Corollary 3.16, Lemma 3.11,
and Lemma 4.2), we obtain\break

\newpage

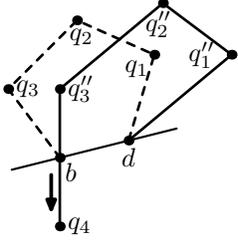
\quad

\vskip4pt

\noindent
$\vcenter{\hbox{\epsfbox{Picture.29}}}$

\leftskip93pt

\vskip-103pt

\noindent
$q_i=q''_i$, $i=1,2,3$, after a few SEs involving $q_1,q_2$ and
$q_2,q_3$. Now, the point $q_3$ is in $\G{\prec}b,q_4{\succ}$. Thus, in
order to execute a given SE of $\varrho_1$ involving $b,q_4$, we can
simply apply a suitable SE of $\varrho$ involving $q_3,q_4$.

In the same manner, we can `cut' $\varrho'$ into $\varrho'_0$ and
$\varrho'_1$ by means of appropriate $b',d'\in\B W$. By induction on
$n$, we can assume that $\varrho_1=\varrho'_1$. It remains to apply
Lemma 4.2
$_\blacksquare$

\medskip

The relation $R(q_5)R(q_4)R(q_3)R(q_2)R(q_1)=1$ suits [Ana1, Conjecture
1.1]. We~strongly believe that this conjecture is valid for the
Poincar\'e disc.

\leftskip0pt

\medskip

Let $\varrho\in\Cal R^\pm H_n$, $\varrho(r_i)=R(q_i)$, $q_i\in\B W$.
Following the proof of Proposition 3.14, associate to $\varrho$ a {\it
standard\/} fundamental polygon $P_\varrho$ for $\varrho H_n$ with
vertices $p_1,p_2,\dots,p_{n-1}$ by taking $p_0:=q_n$ and
$p_i:=R(q_i)p_{i-1}$, $i=1,2,\dots,n-1$. The polygon $P_\varrho$ is
convex and the sum of its interior angles equals $\pm\pi$. In order to
describe $\varrho$, it suffices to mark the vertex $q_n$ and the middle
points $q_1,q_2,\dots,q_{n-1}$ of the edges of $P_\varrho$. Clearly,
$p_{n-1}=p_n=p_0=q_n$. We alter our convention concerning the notation
of the vertices of $P_\varrho$ : the indices of the vertices
$p_1,p_2,\dots,p_{n-1}$ are modulo $n-1$. According to the new
convention, $p_0=p_{n-1}$ and $p_n=p_1$.

\smallskip

We are going to study the group $\Aut H_n$. Fix some discrete subgroup
$H_n\le\Cal L$ and consider the representations
$\varrho\in\Cal R^\pm H_n$ such that $\varrho H_n=H_n$. The group
$\Aut H_n$ acts from the right on these representations. In particular,
every DT can be regarded as an element in $\Aut H_n$ : the automorphism
corresponding to $E_i$ is given by $E_ir_{i-1}=r_i$,
$E_ir_i=r_ir_{i-1}r_i$, and $E_ir_j=r_j$ for $j\notin\{i-1,i\}$.

Denote by $\Aut^+H_n$ the subgroup in $\Aut H_n$ generated by all the
$E_i$'s. In addition, there is an automorphism $J\in\Aut H_n$ given by
$Jr_i:=r_{n-i}$ (cf.~Remark 3.3). Obviously, $J^2=1$.

Define $S\in\Aut H_n$ as $Sr_i:=r_{i+1}$ for all $i$. It is immediate
that $E_i^S=E_{i+1}$. Looking at the polygon~$P_\varrho$,

\vskip3pt

\noindent
\hskip2pt$\vcenter
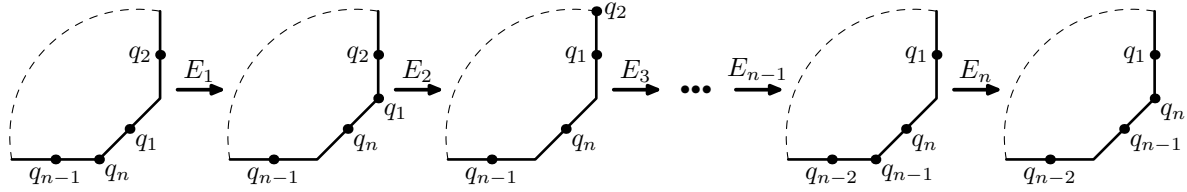
{\hbox{\epsfbox{Picture.30}}}$

\vskip5pt

\noindent
we can see that $S=E_1E_2\dots E_{n-1}\in\Aut^+H_n$. Also, the vertices
$p'_i$ of the standard polygon $P_{\varrho'}$ for the representation
$\varrho':=\varrho SE_n=\varrho E_1E_2\dots E_n$ are given by
$p'_i=p_{i+1}$, where the $p_i$'s stand for the vertices of
$P_\varrho$. Therefore, acting by $\Aut^+H_n$ on the representations,
we can shift the indices both of the vertices and of the marks of the
middle points of the edges of $P_\varrho$.

Denote by $I_h\in\Aut H_n$ the conjugation by $h\in H_n$. Clearly,
$I_h^A=I_{Ah}$ for all $A\in\Aut H_n$. Looking at the polygon
$P_\varrho\cup R(q_1)P_\varrho$,

\vskip6pt

\noindent
\hskip57pt$\vcenter
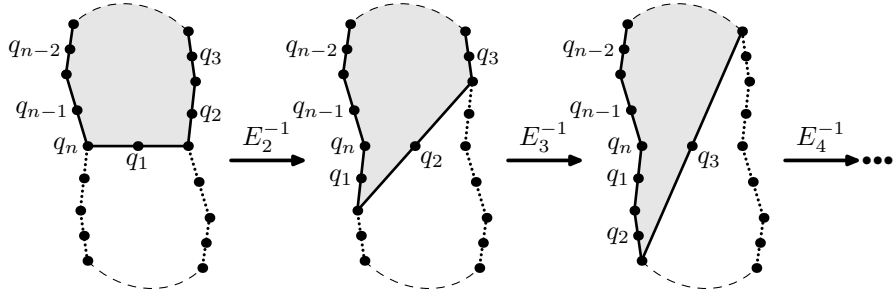
{\hbox{\epsfbox{Picture.32}}}$

\newpage

\noindent
$\vcenter
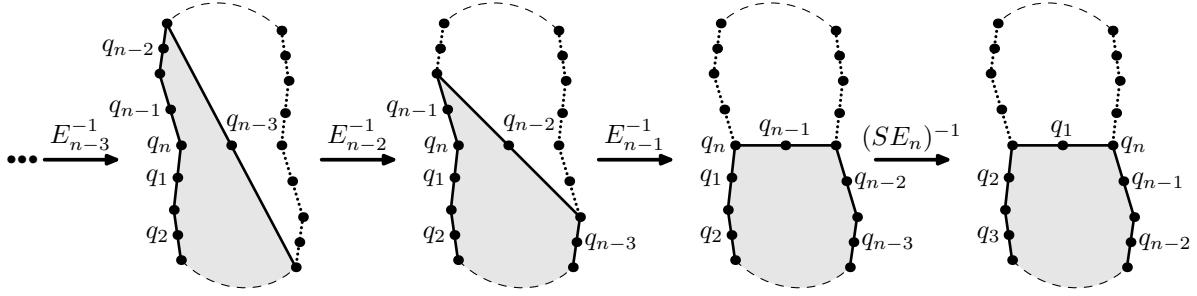
{\hbox{\epsfbox{Picture.33}}}$

\vskip6pt

\noindent
we can see that
$I_{r_1}=E_2^{-1}E_3^{-1}\dots E_{n-1}^{-1}(SE_n)^{-1}$. Hence,
$$I_{r_1}=E_1E_2\dots E_{n-1}E_nE_{n-1}\dots E_3E_2.$$
It follows from $I_{r_i}^S=I_{r_{i+1}}$ that $I_{H_n}\subset\Aut^+H_n$.

\medskip

{\bf4.6.~Theorem.} {\sl The group\/ $\Aut H_n$ is generated by\/ $J$
and by the normal subgroup\/ $\Aut^+H_n$ of index\/~$2$.}

\medskip

{\bf Proof.} Given $\varrho,\varrho'\in\Cal RH_n$ such that
$\varrho H_n=\varrho'H_n=H_n\le\Cal L$, we can assume that
$\varrho,\varrho'\in\Cal R^+H_n$ acting by $J$ if necessary. Hence, the
vertices $p_1,p_2,\dots,p_{n-1}$ and $p'_1,p'_2,\dots,p'_{n-1}$ of the
convex polygons $P:=P_\varrho$ and $P':=P_{\varrho'}$ are listed in
$\partial P$ and in $\partial P'$ in the counterclockwise sense. It
suffices to show that, acting by $\Aut^+H_n$ on both $\varrho$ and
$\varrho'$, we can make them coincide.

\smallskip

Note that the `DT $E$ involving $q_{n-1},q_1$' is expressible in terms
of $E_i$'s : $E=E_1E_nE_1^{-1}$.

\vskip5pt

\noindent
$\vcenter
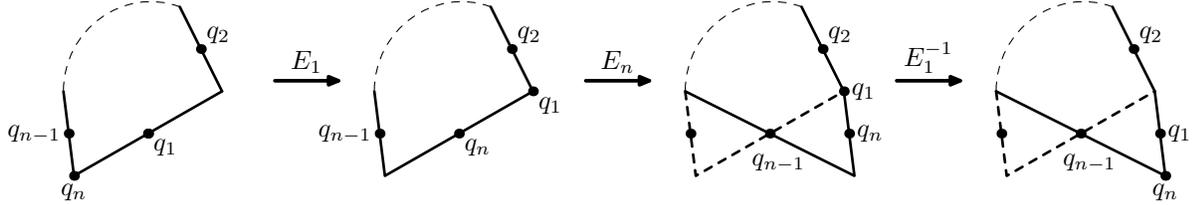
{\hbox{\epsfbox{Picture.31}}}$

\vskip4pt

\noindent
Dealing with the representations $\varrho$ and $\varrho'$ modulo the
action by $\Aut^+H_n$ and taking into account that the automorphism
$SE_n\in\Aut^+H_n$ shifts the indices of the vertices and of the marks
of the middle points of the edges, we can actually think of the
representations as their standard counterclockwise-oriented polygons
$P$ and $P'$, but with unmarked vertices and middle points. As shown
above, we are able to execute any DT that involves the middle points of
adjacent edges of the unmarked polygons, acting by $\Aut^+H_n$. Also,
the inclusion $I_{H_n}\subset\Aut^+H_n$ allows us to change $P$ and
$P'$ by their conjugates.

\smallskip

Let $\varrho(r_i)=R(q_i)$ and $\varrho'(r_i)=R(q'_i)$,
$q_i,q'_i\in\B W$. Every involution $r\in H_n$ is determined by its
fixed point $q$ and induces in $\Sigma$ the hyperelliptic involution
$\iota$. In particular, $\pi q=f_i$ for a suitable $i=1,\dots,n$ (see
the proof of Proposition 4.1). Two involutions $r,r'$ are conjugated in
$H_n$ (equivalently, by an element in $G_n$) if and only if their fixed
points $q,q'$ satisfy the relation $\pi q=\pi q'$. Hence,\break

\vskip-7pt

\noindent
\hskip371pt$\vcenter{\hbox{\epsfbox{Picture.34}}}$

\rightskip88pt

\vskip-79pt

\noindent
the $R(q_i)$'s list all conjugate classes of the involutions in $H_n$.
Obviously, the $R(q'_i)$'s represent different conjugate classes.
Therefore, every $q_i$ is a conjugate of some $q'_j$ and vice versa.

The edge $e_i$ of $P$ has the ends $p_{i-1},p_i$ and the middle point
$q_i$, $i=1,2,\dots,n-1$. Similarly, we introduce the edges $e'_i$ of
$P'$. If $e_i$ is a conjugate of some $e'_j$, we say that $e_i$ and
$e'_j$ are {\it good.} Note that $e_i$ cannot be a conjugate of two
$e'_j$'s at the same time. Let $k$ denote the number of good $e_i$'s.
We proceed by induction on $k\ge0$.

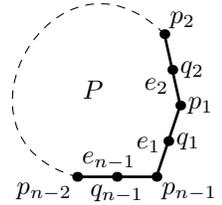
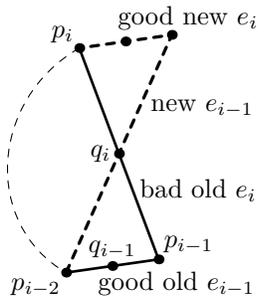
\rightskip0pt

\newpage

\quad

\vskip-7pt

\noindent
$\vcenter{\hbox{\epsfbox{Picture.35}}}$

\leftskip105pt

\vskip-113pt

Suppose that $e_{i-1}$ is a good edge and that $e_i$ is not (the
indices are modulo $n-1$). Apply to $\varrho$ the DT that involves
$q_{i-1},q_i$. This does not alter the edges $e_j$, $j\ne i-1,i$. The
new $e_i$ is a conjugate of the old $e_{i-1}$ and, hence, is good. We
can assume that the new $e_{i-1}$ is bad since, otherwise, we are done
by induction on $k$. So, we are able to permute the types of any two
adjacent edges, one good and the other bad, finally reaching the
situation where the good edges of $P$ (and of $P'$) form a sequence in
$\partial P$ (and~in $\partial P'$). Moreover, we can assume that the
first edge $e$ in the sequence in $\partial P$ and the first edge $e'$
in the sequence in $\partial P'$ are conjugated (both sequences are
read in the counterclockwise sense). By means of DT's, we~can change
$P'$ by any of its conjugates. Also, by means of DT's, we can
shift\break

\vskip-12pt

\leftskip0pt

\noindent
the marks of the vertices and of the middle points in $P$ and in $P'$.
So, we assume that $P$ and $P'$ are on the same side from
$e_1=e=e'=e'_1$. (If $k=0$, we assume only that $p_{n-1}=p'_{n-1}$.)

\vskip8pt

\noindent
\hskip372pt$\vcenter{\hbox{\epsfbox{Picture.36}}}$

\rightskip87pt

\vskip-53pt

The fact that conjugated points in $P$ are necessarily in $\partial P$
and the same fact concerning $P'$ imply that $e_2=e'_2$. In this way,
we can show that $e_i=e'_i$, $i=1,2,\dots,k$. Denote by
$s\subset\partial P$ the segment formed by all the good edges of $P$.
Clearly, $s\subset\partial P'$ is the segment formed by all the good
edges of $P'$. (If $k=0$, we~have $s=p_{n-1}=p'_{n-1}$.)

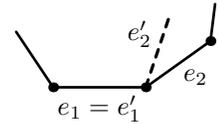
\rightskip0pt

Suppose that $k\ne n$. We will study how the conjugates of bad edges of
$P'$ intersect the polygon $P$. Let $b'\ne b''$ be such edges. Looking
at the tessellation of $\B W$ related to $P'$, we see that

\smallskip

\noindent
{\bf(4.6.1)}~The edges $b'$ and $b''$ can intersect only in points that
are conjugates of $p'_i$'s, i.e., conjugates of $p_{n-1}=p'_{n-1}=q_n$.
Therefore, $b'$ and $b''$ do not intersect in the interior of $P$.

\smallskip

\noindent
{\bf(4.6.2)}~If $b'$ intersects the interior of $P$, it does not
intersect the interior of $s$. Otherwise, $b'$ enters the interior of
$P$ right after its intersection $s\cap b'$ since $P$ is convex. Hence,
it enters the interior of $P'$. A~contradiction.

\smallskip

\noindent
{\bf(4.6.3)}~The edge $b'$ cannot pass through two middle points of
edges of $P$ because the conjugates of middle points of edges of $P$
coincide with those for $P'$.

\smallskip

\noindent
{\bf(4.6.4)}~For every middle point $q_i$ of a bad edge of $P$, there
exists a unique conjugate $b'$ of an edge of~$P'$, necessarily bad,
that passes through $q_i$ and, therefore, through the interior of $P$.

\smallskip

We say that the intersection of $\partial P$ with some conjugate of a
bad edge of $P'$ is {\it proper\/} if this intersection is different
from the vertices of $P$ and from the middle points of the edges of
$P$. It is immediate that the number of proper intersections is the
same in each half of a bad edge of $P$. Let $l$ denote the total number
of proper intersections in $\partial P$. We proceed by induction on
$l$.

\vskip3pt

\noindent
\hskip353pt$\vcenter{\hbox{\epsfbox{Picture.37}}}$

\rightskip105pt

\vskip-82pt

Let $q_i$ be the middle point of a bad edge of $P$ and let $b'$ be a
conjugate of a bad edge of $P'$ that passes through $q_i$ and through
the interior of $P$ according to (4.6.4). By (4.6.2), $b'$ cuts $P$
into two closed parts and $s$ is entirely included in one of them. If
the other part contains a single middle point of an edge of $P$, namely
$q_i$, we arrive at the desired situation to be studied later.
Otherwise, by~(4.6.4), we take a conjugate $b''$ of a bad edge of $P'$
passing through the extra middle point $q_j$ and through the interior
of $P$. Note that $q_j\notin b'$ by (4.6.3).\break

\vskip-12pt

\rightskip0pt

\noindent
By (4.6.1), $b'$ and $b''$ do not intersect in the interior of $P$. Now
we take $b''$ in place of $b'$ and so~on~$\dots$ Finally, we arrive at
the situation (or at the one symmetric to it) where
$b'\cap\partial P=\{q_i,q\}$ and $q\in\G(p_i,q_{i+1})$.

In this situation, we execute the DT $E_{i+1}$. By induction on $k$, we
can assume that the new $e_i$ is bad. We will show that the new $l$ is
strictly less than the old one.

Note that $E_{i+1}$ removes from $P$ the triangle
$\Delta(p_i,q_{i+1},p_{i-1})$ and glue to $P$ the triangle
$\Delta\big(p_{i+1},q_{i+1},\allowmathbreak R(q_{i+1})p_{i-1}\big)$.
Since these triangles are conjugated, it suffices to show that the
number of proper intersections included in $\G(q_{i+1},p_{i-1})$ is
strictly less than that in $\G(p_i,q_{i+1})$. So, we consider only
those parts of conjugates of bad edges of $P'$ that pass via the
interior of $\Delta(p_{i-1},p_i,q_{i+1})$.

\smallskip

The following types and quantities of such parts are possible:

\smallskip

\vskip10pt

\noindent
$\vcenter
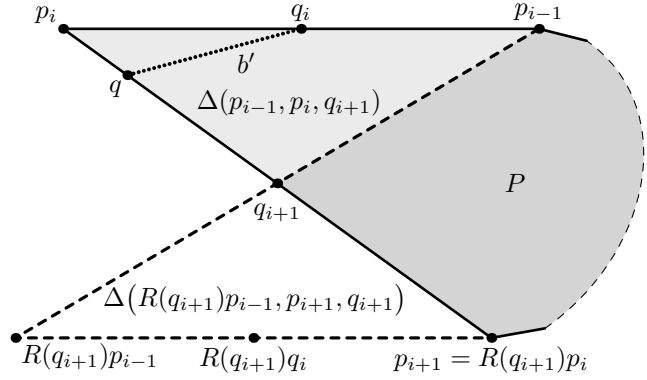
{\hbox{\epsfbox{Picture.38}}}$

\leftskip253pt

\vskip-151pt

\noindent
$\bullet$ $l_1$ parts whose ends are a point in $\G(p_i,q_i)$ and a
point in $\G(p_i,q)$,

\noindent
$\bullet$ 1 part with ends $q_i$ and $q$,

\noindent
$\bullet$ $l_2$ parts whose ends are a point in $\G(q_i,p_{i-1})$ and a
point in $\G(q,q_{i+1})$,

\noindent
$\bullet$ $l_3=0,1$ parts whose ends are a point in $\G(q_i,p_{i-1})$
and $q_{i+1}$,

\noindent
$\bullet$ $l_4$ parts whose ends are a point in $\G(q_i,p_{i-1})$ and a
point in $\G(q_{i+1},p_{i-1})$,

\noindent
$\bullet$ $l_5$ parts whose ends are $p_{i-1}$ and a point in
$\G(q,q_{i+1})$,

\noindent
$\bullet$ $l_6$ parts whose ends are a point in $\G(q,q_{i+1})$ and a
point in $\G(q_{i+1},p_{i-1})$.

\leftskip0pt

\smallskip

Since the number of proper intersections is the same in each half of
$e_i=\G(p_i,p_{i-1})$, we obtain $l_1=l_2+l_3+l_4$. The number of
proper intersections included in $\G(p_i,q_{i+1})$ equals
$l_1+1+l_2+l_5+l_6$. The number of such intersections related to
$\G(q_{i+1},p_{i-1})$ is equal to $l_6+l_4$
$_\blacksquare$

\medskip

A straightforward verification shows that $E_i^J=E_{n+1-i}^{-1}$.
Denote
$$S:=E_1E_2\dots E_{n-1},\qquad\widehat S:=E_{n-1}\dots E_2E_1,\qquad
I:=I_{r_n}.$$
It follows from $I_{r_i}^S=I_{r_{i+1}}$, $E_i^S=E_{i+1}$, and
$I_{r_1}=E_1E_2\dots E_{n-1}E_nE_{n-1}\dots E_3E_2$ that
$I_{r_1}=S{\widehat S}^S$ and $I=I_{r_n}=I_{r_1}^{S^{-1}}=S\widehat S$.
Hence, $S\widehat SS\widehat S=1$. The relations $r_n\dots r_2r_1=1$,
$I_{r_i}^S=I_{r_{i+1}}$, and $S^n=1$ imply the relation
$I^{S^n}\dots I^{S^2}I^S=1$ which can be rewritten as $(S^{-1}I)^n=1$,
i.e., as ${\widehat S}^n=1$. It is immediate that $E_iE_j=E_jE_i$ if
$|i-j|\ge2$. As is easy to see, the relation

\vskip8pt

\noindent
$\vcenter
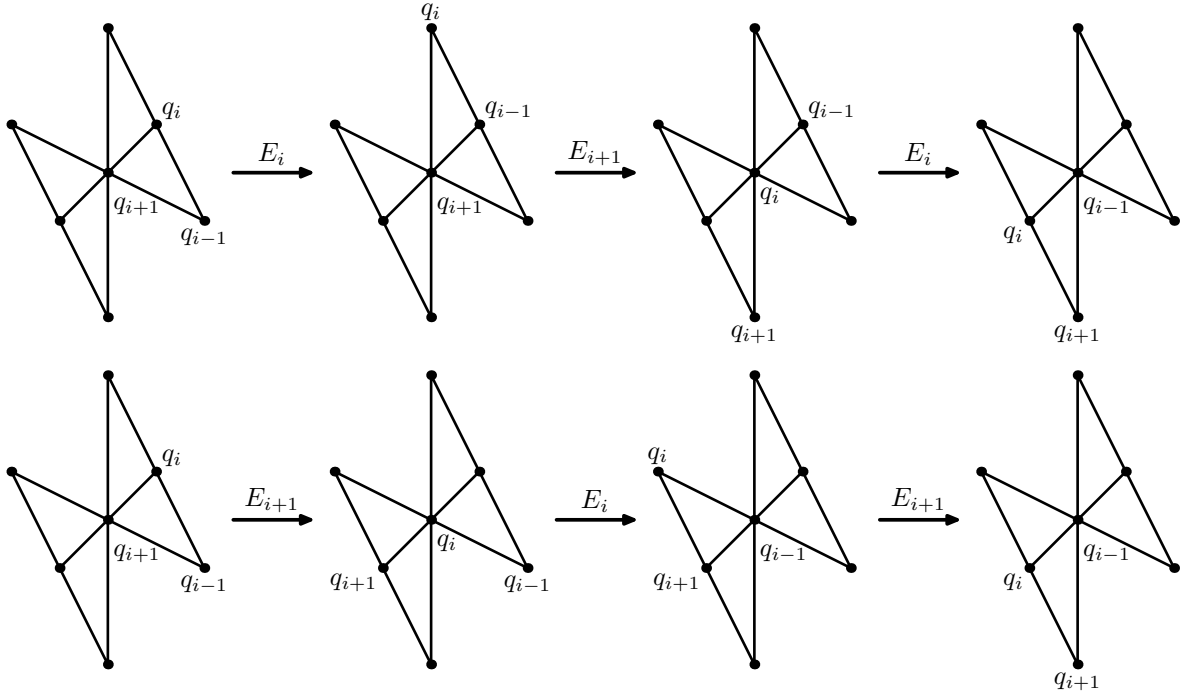
{\hbox{\epsfbox{Picture.39}}}$

\vskip8pt

\noindent
$\vcenter{\hbox{\epsfbox{Picture.40}}}$

\newpage

\noindent
$E_iE_{i+1}E_i=E_{i+1}E_iE_{i+1}$ is valid for all $i$. It is possible
to conclude from [Stu] that the defining relations of $\Aut^+H_n$ are
(the indices are modulo $n$) :
$$S=E_1E_2\dots E_{n-1},\qquad\widehat S=E_{n-1}\dots E_2E_1,\qquad
S^n=1,\qquad\widehat S^n=1,\qquad S\widehat SS\widehat S=1,$$
$$E_i^S=E_{i+1},\qquad E_iE_{i+1}E_i=E_{i+1}E_iE_{i+1},\qquad
E_iE_j=E_jE_i\text{ if }|i-j|\ge2$$
(cf.~[Bir]). The additional defining relations of $\Aut H_n$ are
$E_i^J=E_{n+1-i}^{-1}$ and $J^2=1$.

\bigskip

\centerline{\bf5.~W.~M.~Goldman's Theorem}

\medskip

Let $n\ge6$ be even. Recall that $G_n$ denotes the fully characteristic
torsion-free subgroup of index $2$ in $H_n$ constituted by the words of
even length in the $r_i$'s. By Proposition 4.1, $G_n$ is the
fundamental group of a closed orientable Riemann surface of genus
$\frac n2-1$. In this section, we will prove the

\medskip

{\bf5.1.~Theorem {\rm[Gol1, Corollary C]}.} {\sl Let\/
$\varrho:G_n\to\Cal L$ be a representation. Then\/
$\varrho\in\Cal RG_n$ if and only if\/ $\Area\varrho=\pm2(n-4)\pi$.}

\medskip

We are going to explore the ideas developed in the hyperelliptic case.
A given representation $\varrho:G_n\to\Cal L$ defines an action of
$G_n$ on $\overline\B W$. We write $gp$ instead of $\varrho(g)p$ for
all $g\in G_n$ and $p\in\overline\B W$. Working in terms of the
$r_i$'s, we are allowed to apply $\varrho$ to any expression of even
length in $r_i$'s. Hence, the expression $r_ir_jp$ makes sense, whereas
$r_ip$ does not.

\vskip5pt

\noindent
\hskip325pt$\vcenter{\hbox{\epsfbox{Picture.41}}}$

\rightskip133pt

\vskip-103pt

We will deal with a `fundamental polygon' $Q$ for $\varrho G_n$ that
mimics the duplicated fundamental polygon $P_n$ for the hyperelliptic
case, namely, $Q:=P_n\cup\varrho(r_n)P_n$ (see the last picture in the
proof of Theorem 3.15). In the hyperelliptic case, the polygon $P_n$ is
generated by the choice of $p=p_n\in\B W$ because it has a single cycle
of vertices. The point $p_{n-1}\in\B W$ is given by
$p_{n-1}=\varrho(r_n)p_n$. Since, in the nonhyperelliptic case, we have
no reflection $\varrho(r_n)$ available and the polygon $Q$ should have
two cycles of vertices, we choose two points $p,q\in\overline\B W$ that
are intended to respectively play the roles of $p_n,p_{n-1}$. In this
way, for suitable $w_i\in G_n$,\break

\vskip-12pt

\rightskip0pt

\noindent
the even vertices of the polygon $Q$ have the form $w_{2j}p$ and the
odd ones, the form $w_{2j+1}q$.

\smallskip

The proof of Theorem 5.1 is `almost' the same as that of Theorem 3.15.
We simply adapt the arguments of the latter to the nonhyperelliptic
case by avoiding the use of the elements from $H_n\setminus G_n$. For
instance, Corollary 5.8, Remark 5.9, Remark 5.10, Lemma 5.12, and
Lemma 5.13 that we prove below are analogs of the following
hyperelliptic assertions: Lemma 3.2, Remark 3.3, Remark 3.4, Lemma~3.9,
and Proposition 3.13.

\medskip

{\bf5.2.~Notation.} Denote by $S$, $I$, and $J$ the automorphisms of
$H_n$ given by the rules $Sr_i=r_{i+1}$, $I:h\mapsto h^{r_n}$, and
$Jr_i:=r_{n-i}$. The same symbols denote the induced automorphisms of
$G_n$. For~$0\le i\le n-1$, denote $v_i:=r_i\dots r_2r_1$ and regard
the indices of the $v_i$'s modulo $n$. So, $v_0=v_n=1$.
For~$0\le i\le n-2$, introduce
$$w_i:=v_i\text{ if }i\text{ is even},\qquad w_i:=v_ir_n\text{ if
}i\text{ is odd},\qquad w_{i+n-1}:=I(w_i)$$
and regard the indices of the $w_i$'s modulo $2n-2$. Clearly,
$w_0=w_{n-1}=w_{2n-2}=1$. Note that $w_{i+n-1}=I(w_i)$ for all $i$. As
is easy to see, the formula $w_i=v_ir_n$ works for all odd $i$ such
that $1\le i\le n-1$.

\medskip

The elementary properties of the $w_i$'s that we use in what follows
are gathered in the

\newpage

{\bf5.3.~Lemma.}

{\sl{\rm(1)}\/ $w_{i+n-1}^{-1}w_{i+n}=w_{i+1}^{-1}w_i$ for all\/ $i$.

{\rm(2)}\/ $J(w_i)=w_{n-1-i}$ for all\/ $i$.

{\rm(3)}\/ $S(w_i)w_1=w_{i+1}$ for all even\/ $i$ such that\/
$0\le i\le n-2$.

{\rm(4)}\/ $S(w_i)=w_{i+1}$ for all odd\/ $i$ such that\/
$1\le i\le n-3$.

{\rm(5)}\/ $S(w_i)=w_1w_{i+1}$ for all odd\/ $i$ such that\/
$n-1\le i\le2n-3$.

{\rm(6)}\/ $S(w_i)w_1=w_1w_{i+1}$ for all even\/ $i$ such that\/
$n\le i\le2n-4$.

{\rm(7)}\/ $r_nr_iw_{i-1}=w_{i+n-1}$ and $r_nr_iw_i=w_{i+n-2}$ for
all\/ $1\le i\le n-1$.

{\rm(8)}\/ $r_nr_{i+1}w_i=w_{i+n}$, $r_nr_{i+1}w_{i+1}=w_{i+n-1}$,
$r_{i+1}r_nw_{i+n-1}=w_{i+1}$, and $r_{i+1}r_nw_{i+n}=w_i$ for all\/
$2\le i\le n$.}

\medskip

{\bf Proof.} (1)~Let $0\le i\le n-2$. If $i$ is even, we have
$$w_{i+n-1}^{-1}w_{i+n}=(r_nw_ir_n)^{-1}r_nw_{i+1}r_n=
r_nv_i^{-1}v_{i+1}=r_nv_i^{-1}r_{i+1}v_i=w_{i+1}^{-1}w_i.$$
If $i$ is odd, we have
$$w_{i+n-1}^{-1}w_{i+n}=(r_nw_ir_n)^{-1}r_nw_{i+1}r_n=
v_i^{-1}v_{i+1}r_n=v_i^{-1}r_{i+1}v_ir_n=w_{i+1}^{-1}w_i.$$
For $n-1\le i\le2n-3$, the fact follows by taking inverses in the
equalities that are already established for $0\le i\le n-2$.

(2)~Let $0\le i\le n-2$. It follows from the relation
$r_n\dots r_2r_1=1$ that
$$J(w_i)=J(v_i)=r_{n-i}\dots r_{n-2} r_{n-1}=v_{n-1-i}r_n=w_{n-1-i}$$
if $i$ is even and that
$$J(w_i)=J(v_ir_n)=r_{n-i}\dots r_{n-2} r_{n-1}r_n=v_{n-1-i}=
w_{n-1-i}$$
if $i$ is odd. Now, for $n-1\le i\le2n-3$, we obtain
$$J(w_i)=J(r_nw_{i-n+1}r_n)=r_nw_{n-1-i+n-1}r_n=r_nw_{2n-2-i}r_n=
w_{3n-3-i}=w_{n-1-i}.$$

(3)~The case of $i=0$ is immediate. For $2\le i\le n-2$, we have
$S(w_i)w_1=S(v_i)r_1r_n=v_{i+1}r_n=w_{i+1}$.

(4)~$S(w_i)=S(v_ir_n)=v_{i+1}=w_{i+1}$.

(5)~$S(w_i)=S(r_nw_{i-n+1}r_n)=S(r_nv_{i-n+1}r_n)=r_1v_{i-n+2}=
r_1w_{i-n+2}r_n=r_1r_nw_{i+1}=w_1w_{i+1}$.

(6)~$S(w_i)w_1=S(r_nw_{i-n+1}r_n)w_1=S(r_nv_{i-n+1})r_1r_n=
r_1v_{i-n+2}r_n=r_1w_{i-n+2}r_n=w_1w_{i+1}$.

(7)~As is easy to see, $r_iw_{i-1}=w_ir_n$ and $r_iw_i=w_{i-1}r_n$ for
all $1\le i\le n-2$. Therefore,
$r_nr_iw_{i-1}=r_nw_ir_n=I(w_i)=w_{i+n-1}$ and
$r_nr_iw_i=r_nw_{i-1}r_n=I(w_{i-1})=w_{i+n-2}$. For $i=n-1$, we have
$r_nr_{n-1}w_{n-2}=1=w_{2n-2}$ and
$r_nr_{n-1}w_{n-1}=r_nr_{n-1}=r_nv_{n-2}r_n=I(w_{n-2})=w_{2n-3}$ since
$r_nr_{n-1}\dots r_2r_1=1$, $r_{n-1}=r_{n-2}\dots r_2r_1r_n$, and
$w_{n-1}=1$.

(8)~The first two equalities are in fact shown in (7). The last two
equalities follow immediately from the first two
$_\blacksquare$

\medskip

Given $p,q\in\overline\B W$, define
$${\Area}_n(p,q;\varrho):=\Area(w_0p,w_1q,\dots,w_{n-2}p,w_{n-1}q,
\dots,w_{2n-4}p,w_{2n-3}q),$$
$${\Area}_{i+1}(p,q;\varrho):={\Area}_i(p,q;\varrho S).$$

\medskip

{\bf5.4.~Remark.} The relation $w_{i+n-1}=r_nw_ir_n$ valid for all $i$
implies $\Area_n(p,q;\varrho)=\Area_n(q,p;\varrho I)$.

\medskip

{\bf5.5.~Lemma.} {\sl$\Area_n(p,q;\varrho)=\Area_1(w_1q,p;\varrho)$.}

\medskip

{\bf Proof.} By definition,
$${\Area}_1(w_1q,p;\varrho)=\Area\big(S(w_0)w_1q,S(w_1)p,\dots,
S(w_{n-2})w_1q,S(w_{n-1})p,\dots,S(w_{2n-4})w_1q,S(w_{2n-3})p\big).$$
By Lemma 5.3 (3--6),
$${\Area}_1(w_1q,p;\varrho)=\Area(w_1q,w_2p,\dots,w_{n-1}q,w_1w_np,
\dots,w_1w_{2n-3}q,w_1w_{2n-2}p).$$
Taking into account that $w_1w_n=w_1r_nw_1r_n=1$ and that
$w_0=w_{n-1}=w_{2n-2}=1$, by Remark 2.3, we obtain
$${\Area}_1(w_1q,p;\varrho)=\Area(w_1q,w_2p,\dots,w_{n-1}q,w_1w_np)+
\Area(w_1w_np,\dots,w_1w_{2n-3}q,w_1w_{2n-2}p,w_1q)=$$
$$=\Area(w_0p,w_1q,w_2p,\dots,w_{n-1}q)+
\Area(w_np,\dots,w_{2n-2}p,q)=$$
$$=\Area(w_0p,w_1q,\dots,w_{n-2}p,w_{n-1}q)+
\Area(w_{n-1}q,w_np,w_{n+1}q,\dots,w_{2n-3}q,w_0p)=$$
$$=\Area(w_0p,w_1q,\dots,w_{n-2}p,w_{n-1}q,\dots,w_{2n-3}q)=
{\Area}_n(p,q;\varrho)\ _\blacksquare$$

\medskip

{\bf5.6.~Lemma.} {\sl$\Area_n(p,q;\varrho)$ is independent of the
choice of\/ $p$ and\/ $q$.}

\medskip

{\bf Proof.} We will show the independence of $q$. (The independence of
$p$ can be shown in a similar way.) Taking $c=p$ in (2.2), we obtain
$${\Area}_n(p,q;\varrho)=\sum_{\text{even
}i}\Area\Delta(p,w_ip,w_{i+1}q)+\sum_{\text{odd
}i}\Area\Delta(p,w_iq,w_{i+1}p)=$$
$$=\sum_{\text{even
}i}\Area\Delta(q,w_{i+1}^{-1}p,w_{i+1}^{-1}w_ip)+\sum_{\text{odd
}i}\Area\Delta(q,w_i^{-1}w_{i+1}p,w_i^{-1}p).\leqno{\bold{(5.7)}}$$
Let us show that (5.7) is the area (calculated with respect to the
centre $q$) related to some closed piecewise geodesic path $C$
independent of the choice of $q$. Denote by
$\overset i\to\longrightarrow$ the side opposite to the vertex $q$ of
the $i$th triangle involved in (5.7). This side is oriented with
respect to the orientation of the $i$th triangle. The consecutive
vertices of $C$ are described by the following list:
$$w_1^{-1}p\overset0\to\longrightarrow
w_1^{-1}w_0p=w_{n-1}^{-1}w_np\overset n-1\to\longrightarrow
w_{n-1}^{-1}p\overset n-2\to\longrightarrow
w_{n-1}^{-1}w_{n-2}p=w_{2n-3}^{-1}w_0p\overset2n-3\to\longrightarrow
w_{2n-3}^{-1}p\overset2n-4\to\longrightarrow\dots$$
$$\dots\overset2j\to\longrightarrow w_{2j+1}^{-1}w_{2j}p=
w_{2j+n-1}^{-1}w_{2j+n}p\overset2j+n-1\to\longrightarrow
w_{2j+n-1}^{-1}p\overset2j+n-2\to\longrightarrow
w_{2j+n-1}^{-1}w_{2j+n-2}p=
w_{2j-1}^{-1}w_{2j}p\overset2j-1\to\longrightarrow\dots$$
$$\dots\overset3\to\longrightarrow w_3^{-1}p\overset2\to\longrightarrow
w_3^{-1}w_2p=w_{n+1}^{-1}w_{n+2}p\overset n+1\to\longrightarrow
w_{n+1}^{-1}p\overset n\to\longrightarrow
w_{n+1}^{-1}w_np=w_1^{-1}w_2p\overset1\to\longrightarrow w_1^{-1}p,$$
where the equalities are provided by Lemma 5.3 (1). In this list, the
mentioned sides of even triangles appear in the order
$$\overset0\to\longrightarrow\qquad\overset
n-2\to\longrightarrow\qquad\dots\qquad\overset2j\to\longrightarrow
\qquad\overset2j+(n-2)\to\longrightarrow\qquad\dots\qquad\overset
n\to\longrightarrow$$
and the mentioned sides of odd ones, in the order
$$\overset n-1\to\longrightarrow\qquad\overset
n-1+(n-2)\to\longrightarrow\qquad\dots\qquad\overset2j+1\to
\longrightarrow\qquad\overset2j+1+(n-2)\to\longrightarrow\qquad\dots
\qquad\overset1\to\longrightarrow.$$
Since $n-2$ and $n-1$ are coprime, every side appears exactly once in
the list
$_\blacksquare$

\medskip

{\bf5.8.~Corollary.} {\sl$\Area_i(p,q;\varrho)$ does not depend on the
choice of\/ $p$, $q$, and\/ $i$
$_\blacksquare$}

\medskip

{\bf5.9.~Remark.} By Lemma 5.3 (2), $\Area\varrho J=-\Area\varrho$.

\medskip

In the sequel, we assume without loss of generality that
$\Area\varrho\ge0$.

\medskip

\noindent
\hskip294pt$\vcenter{\hbox{\epsfbox{Picture.42}}}$

\vskip-114pt

\rightskip165pt

{\bf5.10.~Remark.} Take a fixed point $c=p=q\in\overline\B W$ of
$\varrho(w_1)$. It follows from $w_{n-1}=w_0=1$ and $w_n=w_1^{-1}$ that
$w_1q=w_{n-1}q=w_np=w_0p=c$. Therefore,

\vskip10pt

\hskip33pt$\Area\Delta(c,w_0p,w_1q)=\Area\Delta(c,w_1q,w_2p)=$

\vskip10pt

\hskip12pt$=\Area\Delta(c,w_{n-2}p,w_{n-1}q)=
\Area\Delta(c,w_{n-1}q,w_np)=$

\vskip10pt

\hskip11pt$=\Area\Delta(c,w_np,w_{n+1}q)=\Area\Delta(c,w_{2n-3}q,w_0p)=
0$.

\vskip6pt

\noindent
Hence, $\Area\varrho\le2(n-4)\pi$. When $\Area\varrho=2(n-4)\pi$, we
say that $\Area\varrho$ is {\it maximal.} In this case, $p\in\S W$ and
the cycles

\rightskip0pt

$$p,w_2p,w_3p,\dots,w_{n-2}p,\qquad
p,w_{n+1}p,w_{n+2}p,\dots,w_{2n-3}p$$
are positive.

\medskip

\noindent
\hskip370pt$\vcenter{\hbox{\epsfbox{Picture.43}}}$

\rightskip90pt

\vskip-73pt

{\bf5.11.~Remark.} Let $p_1,p_2,q_2,q_1\in\S W$ be a positive cycle and
suppose that some isometry $h\in\Cal L$ maps $p_i$ to $q_i$, $i=1,2$.
Then $h$ is hyperbolic and the cycle $p_1,s,p_2,q_2,t,q_1$ is positive,
where $s\in\S W$ and $t\in\S W$ stand for the repeller and for the
attractor of $h$.

\medskip

{\bf5.12.~Lemma.} {\sl Let\/ $\varrho:G_n\to\Cal L$ be a representation
with maximal\/ $\Area\varrho$ and let\/ $d\in\S W$ be a fixed point
of\/ $\varrho(w_1)$. Then the cycles}

\rightskip0pt

{\sl
$$d,w_2d,r_3r_1d,w_3d,r_4r_1d,w_4d,\dots,r_{n-3}r_1d,w_{n-3}d,
r_{n-2}r_1d,w_{n-2}d$$
and
$$d,w_{n+1}d,r_nr_3d,w_{n+2}d,r_nr_4d,w_{n+3}d,\dots,r_nr_{n-3}d,
w_{2n-4}d,r_nr_{n-2}d,w_{2n-3}d$$
are positive.}

\medskip

{\bf Proof.} The cycles $d,w_id,w_{i+1}d$ and $d,w_{i+n-1}d,w_{i+n}d$
are positive for all $2\le i\le n-3$ by Remark ~5.10. Hence, by Lemma
5.3 (8), the cycles $r_nr_{i+1}d,w_{i+n}d,w_{i+n-1}d$ and
$r_{i+1}r_nd,w_{i+1}d,w_id$ are positive. In other words, the cycles
$w_{n+i-1}d,r_nr_{i+1}d,w_{n+i}d$ and $w_id,r_{i+1}r_nd,w_{i+1}d$ are
positive. Since $d$ is a fixed point of $w_1=r_1r_n$, we have
$r_{i+1}r_1d=r_{i+1}r_1r_1r_nd=r_{i+1}r_nd$. Therefore, the cycles
$w_id,r_{i+1}r_1d,w_{i+1}d$ and $w_{n+i-1}d,r_nr_{i+1}d,w_{n+i}d$ are
positive for all $2\le i\le n-3$. By Remarks 5.10 and 3.7, the cycles
in Lemma 5.12 are positive
$_\blacksquare$

\medskip

{\bf5.13.~Lemma.} {\sl In the situation of Lemma\/ {\rm5.12}, the
isometry\/ $h_i:=\varrho(r_ir_{i-1})$ is hyperbolic for all\/ $i$
{\rm(}the indices are modulo\/ $n$\/{\rm)}. Denote by\/ $s_{i-1}$ and\/
$t_i$ the repeller and the attractor of\/ $h_i$. Then, for every\/
$d\in\{s_n,t_1\}$, the cycle\/
$$t_1,s_2,w_2d,s_3,t_3,w_3d,s_4,t_4,w_4d,\dots,s_{n-2},t_{n-2},
w_{n-2}d,t_{n-1},s_n$$
is positive.}

\newpage

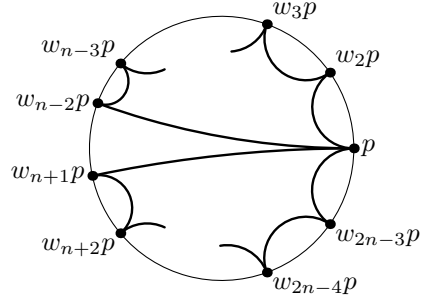
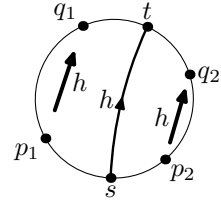
\quad

\vskip-5pt

\noindent
\hskip286pt$\vcenter{\hbox{\epsfbox{Picture.44}}}$

\rightskip172pt

\vskip-142pt

{\bf Proof.} The cycle $w_{i-2}d,r_{i-1}r_1d,r_ir_1d,w_id$ is positive
for all $3\le i\le n-1$. Indeed, for $4\le i\le n-2$, this follows
straightforwardly from Lemma 5.12. For $i=3$, the cycle has the form
$d,w_2d,r_3r_1d,w_3d$ because $w_1d=d$ and $r_2r_1=w_2$. It is positive
by Lemma 5.12. The relation $r_nr_{n-1}\dots r_2r_1=1$ implies
$w_{n-2}=v_{n-2}=r_{n-1}r_n$. From $d=w_1d$ and $w_1=r_1r_n$, we obtain
$r_{n-1}r_1d=r_{n-1}r_1r_1r_nd=r_{n-1}r_nd=w_{n-2}d$. Taking
$w_{n-1}=1$ into account, we can see that, for $i=n-1$, the cycle has
the form $w_{n-3}d,r_{n-2}r_1d,w_{n-2}d,d$. By Lemma 5.12, it is
positive.

The isometry $h_i$ maps $r_{i-1}r_1d$ to $r_ir_1d$ and $w_{i-2}d$ to
$w_id$ for\break

\vskip-12pt

\rightskip0pt

\noindent
all $3\le i\le n-1$. By Remark 5.11, $h_i$ is hyperbolic and the cycle
$$w_{i-2}d,s_{i-1},r_{i-1}r_1d,r_ir_1d,t_i,w_id,\qquad3\le i\le
n-1,\leqno{\bold{(5.14)}}$$
is positive.

The cycle $r_{i-1}r_1d,w_{i-1}d,r_ir_1d$ is positive for all
$4\le i\le n-2$ by Lemma 5.12. We can combine this cycle and the cycle
(5.14) by Remark 3.7 and obtain the positive cycle
$w_{i-2}d,s_{i-1},r_{i-1}r_1d,w_{i-1}d,r_ir_1d,\allowmathbreak
t_i,w_id$
for all $4\le i\le n-2$. The first and the second parts of this cycle
provide the positive cycles
$$w_{i-1}d,s_i,r_ir_1d,w_id,\qquad3\le i\le n-3,\leqno{\bold{(5.15)}}$$
$$w_{i-1}d,r_ir_1d,t_i,w_id,\qquad4\le i\le n-2.\leqno{\bold{(5.16)}}$$

Combining the cycles (5.15) and (5.16) by Remark 3.7, we get the
positive cycle
$$w_{i-1}d,s_i,t_i,w_id,\qquad4\le i\le n-3.\leqno{\bold{(5.17)}}$$

Taking into account that $w_1d=d$ and $r_2r_1=w_2$, we can see that
$d,s_2,w_2d,r_3r_1d,t_3,w_3d$ and $w_2d,s_3,r_3r_1d,w_3d$ are the
cycles (5.14) and (5.15) with $i=3$. Combining these cycles by Remark
3.7 and excluding the term $r_3r_1d$, we arrive at the positive cycle
$$d,s_2,w_2d,s_3,t_3,w_3d.\leqno{\bold{(5.18)}}$$

As was shown above, $r_{n-1}r_1d=w_{n-2}d$. Taking the cycle (5.16)
with $i=n-2$ and the cycle (5.14) with $i=n-1$, we obtain the positive
cycles $w_{n-3}d,r_{n-2}r_1d,t_{n-2},w_{n-2}d$ and
$w_{n-3}d,s_{n-2},r_{n-2}r_1d,\allowmathbreak w_{n-2}d,t_{n-1},d$ since
$w_{n-1}=1$. Combining these cycles by Remark 3.7 and excluding the
term $r_{n-2}r_1d$, we arrive at the positive cycle
$$w_{n-3}d,s_{n-2},t_{n-2},w_{n-2}d,t_{n-1},d.\leqno{\bold{(5.19)}}$$

The cycle $d,w_2d,w_3d,\dots,w_{n-2}d$ is positive by Lemma 5.12.
Combining this cycle with the cycles (5.18), (5.17) for all $i$, and
(5.19), we get the positive cycle
$$d,s_2,w_2d,s_3,t_3,w_3d,s_4,t_4,w_4d,\dots,w_{n-3}d,s_{n-2},t_{n-2},
w_{n-2}d,t_{n-1}.\leqno{\bold{(5.20)}}$$

Shifting the indices, i.e., applying the results already obtained to
the representations $\varrho S^j$, we conclude that $h_i$ is hyperbolic
for all $i$. So, the points $s_n,t_1,s_1,t_2,s_{n-1},t_n$ make sense.

Since the cycle (5.20) is positive for $d=t_1$, the cycle
$t_1,s_2,t_3,s_4$ is positive. Shifting the indices, we conclude that
the cycle $t_{n-1},s_n,t_1,s_2$ is positive. Combining the positive
cycles $t_1,s_2,t_3,s_4$, (5.20), and $t_{n-1},s_n,t_1,s_2$, we arrive
at the positive cycle in Lemma 5.13
$_\blacksquare$

\medskip

{\bf5.21.~Proposition.} {\sl Let\/ $\varrho:G_n\to\Cal L$ be a
representation with maximal\/ $\Area\varrho$. Then the isometries\/
$h_i:=\varrho(r_ir_{i-1})$ and\/ $h'_i:=\varrho(r_nr_ir_{i-1}r_n)$ are
hyperbolic for all\/ $i$ {\rm(}the indices are modulo\/ $n$\/{\rm)}.
Denote by\/ $s_{i-1},s'_{i-1}$ and\/ $t_i,t'_i$ the repellers and the
attractors of\/ $h_i,h_i'$, respectively. Then\/ $s_n=t'_1$,
$t_1=s'_n$, and, for every\/ $d\in\{s_n,t_1\}$, the cycle
$$t_1,s_2,w_2d,s_3,t_3,w_3d,s_4,t_4,w_4d,\dots,s_{n-2},t_{n-2},
w_{n-2}d,t_{n-1},s_n,$$
$$s'_2,w_{n+1}d,s'_3,t'_3,w_{n+2}d,s'_4,t'_4,w_{n+3}d,\dots,s'_{n-2},
t'_{n-2},w_{2n-3}d,t'_{n-1}$$
is positive.}

\medskip

{\bf Proof.} By Lemma 5.13, the isometries $h_i$'s are hyperbolic and
the cycle
$$t_1,s_2,w_2d,s_3,t_3,w_3d,s_4,t_4,w_4d,\dots,s_{n-2},t_{n-2},
w_{n-2}d,t_{n-1},s_n$$
is positive for every $d\in\{s_n,t_1\}$. By Remark 5.4 and Corollary
5.8, $\Area\varrho I=2(n-4)\pi$. By Lemma~5.13 applied to the
representation $\varrho I$, the isometries $h'_i$'s are hyperbolic and
the cycle
$$t'_1,s'_2,w_{n+1}d,s'_3,t'_3,w_{n+2}d,s'_4,t'_4,w_{n+3}d,\dots,
s'_{n-2},t'_{n-2},w_{2n-3}d,t'_{n-1},s'_n$$
is positive for every $d\in\{s'_n,t'_1\}$ since $w_{i+n-1}=I(w_i)$ for
all $i$. It remains to observe that ${h'_1}^{-1}=h_1$ and to combine
the above positive cycles
$_\blacksquare$

\medskip

{\bf Proof of Theorem 5.1.} Let us show that $\Area\varrho=2(n-4)\pi$
implies $\varrho\in\Cal RG_n$.

Denote $\G:=\G[s_n,t_1]$, $\G_i:=w_i\G$, and $\G'_i:=w_{i+n-1}\G$ for
all $2\le i\le n-2$. By Proposition 5.21, $\G$ is the axis of
$h_1={h'_1}^{-1}$. Hence, the vertices of $\G_i$ and of $\G'_i$ are
respectively of the form $w_id$ and $w_{i+n-1}d$, where
$d\in\{s_n,t_1\}$.

Take $p,q\in\G\cap\B W$ such that $p=h_1q$ and denote by $Q$ the
polygon with the successive vertices
$$w_1q,w_2p,w_3q,\dots,w_{n-3}q,w_{n-2}p,w_np,w_{n+1}q,w_{n+2}p,\dots,
w_{2n-4}p,w_{2n-3}q$$
and the successive edges
$e_2,e_3,\dots,e_{n-1},e'_2,e'_3,\dots,e'_{n-1}$ such that
$$e_i:=\G[w_{i-1}q,w_ip],\qquad
e'_i:=\G[w_{i+n-2}p,w_{i+n-1}q]\qquad\text{for even }i,\qquad2\le i\le
n-2,\leqno{\bold{(5.22)}}$$
$$e_i:=\G[w_{i-1}p,w_iq],\qquad
e'_i:=\G[w_{i+n-2}q,w_{i+n-1}p]\qquad\text{for odd }i,\qquad3\le i\le
n-1.\leqno{\bold{(5.23)}}$$
(Note that $w_{n-1}q=w_np$ and $w_{2n-2}p=w_1q$ since
$w_{n-1}=w_{2n-2}=1$ and $w_n=w_1^{-1}$.)

\smallskip

We claim that $Q$ is a fundamental polygon for the group $\varrho G_n$.
Obviously, $w_ip,w_iq\in\G_i$ and $w_{i+n-1}p,w_{i+n-1}q\in\G'_i$ for
all $2\le i\le n-2$. Also, $w_np,w_1q\in\G$ since $w_n=w_1^{-1}$,
$h_1=\varrho(w_1)$, $p=h_1q$, and $p,q\in\G$. Let $d\in\{s_n,t_1\}$.
Then the cycle in Proposition 5.21 is positive. This implies that $\G$,
the $\G_i$'s, and the $\G'_j$'s are all disjoint. Therefore, the edges
$e_i$'s and $e'_i$'s are not degenerated and, thus, generate complete
geodesics $\Gamma_i$ and $\Gamma'_i$.

\smallskip

Define the arcs
$$A:=\{b\in\S W\mid\text{\rm the cycle }t_1,b,s_n\text{ \rm is
positive}\},\qquad A':=\{b\in\S W\mid \text{\rm the cycle
}s_n,b,t_1\text{ \rm is positive}\}.$$
Let $A_i\subset A$ and $A'_i\subset A'$ be the arcs with the same ends
as $\G_i$ and $\G'_i$, respectively. The arcs\break

\vskip-2pt

\noindent
$\vcenter{\hbox{\epsfbox{Picture.45}}}$

\vskip7pt

\noindent
$A',A_2,A_3,\dots,A_{n-2}$ are disjoint because the cycle in
Proposition 5.21 is positive. It is easy to see that the vertices of
$\Gamma_i$ belong to $A_{i-1}$ and $A_i$ for all $3\le i\le n-2$, that
the vertices of $\Gamma_2$ belong to $A'$ and $A_2$, and that the
vertices of $\Gamma_{n-1}$ belong to $A_{n-2}$ and $A'$. The only
intersections between $\Gamma_i$'s are the known intersections between
$\Gamma_{i-1}$ and $\Gamma_i$, $3\le i\le n-1$, and a possible
intersection between $\Gamma_2$ and $\Gamma_{n-1}$. Nevertheless, the
edges $e_2$ and $e_{n-1}$ do not intersect. Indeed, it follows from
Proposition 5.21 that the cycle $t_1,w_2d,s_3,w_{n-2}d,s_n$ is positive
for every $d\in\{s_n,t_1\}$. Since $s_n$ and $t_1$ are the repeller and
the attractor of $h_1=\varrho(w_1)$, $p=h_1q$, $w_np=q$, and $w_1q=p$,
the edges $e_2$ and $e_{n-1}$ cannot intersect. Consequently, the edges
$e_2,e_3,\dots,e_{n-1}$ intersect in the `prescribed' way and are on
the side of the normal vector to $\G$. For similar reasons, the edges
$e'_2,e'_3,\dots,e'_{n-1}$ intersect in the `prescribed' way and are on
the opposite side of the normal vector to $\G$. In other words, $Q$ is
simple.

The polygon $Q$ has $2(n-2)$ vertices and
$\Area\varrho=\Area Q=2(n-4)\pi$. Therefore, the sum of the interior
angles of $Q$ equals $2(n-3)\pi-\Area P=2\pi$. The isometry
$\gamma_i:=\varrho(r_nr_i)$ maps the edge $e_i$ onto the edge $e'_i$
for all $2\le i\le n-1$. This follows from (5.22--23) and from Lemma
5.3 (7).

As is easy to see, the identifications by the $\gamma_i$'s produce the
only cycle of vertices. By Poincar\'e's Polyhedron Theorem, $Q$ is a
fundamental polygon for the group generated by the $\gamma_i$'s and
$$\gamma_{n-1}\dots\gamma_4^{-1}\gamma_3\gamma_2^{-1}\gamma_{n-1}^{-1}
\dots\gamma_4\gamma_3^{-1}\gamma_2=1$$
is a unique defining relation of this group. In other words, $\varrho$
is an isomorphism and, thus, $\varrho\in\Cal RG_n$.

For the converse, we simply repeat the arguments presented at the end
of the proof of Theorem 3.15~$_\blacksquare$

\medskip

{\bf5.24.~Remark.} It is easy to verify that the group $G_n$ admits the
generators $g_{i(i-1)}$ (the indices are modulo $n$) subject to the
defining relations
$$g_{n(n-1)}g_{(n-1)(n-2)}g_{(n-2)(n-3)}\dots
g_{32}g_{21}g_{1n}=1,$$
$$g_{n(n-1)}g_{(n-2)(n-3)}\dots g_{43}g_{21}=1,\qquad
g_{(n-1)(n-2)}g_{(n-3)(n-4)}\dots g_{32}g_{1n}=1.$$
(In terms of $H_n$, $g_{i(i-1)}:=r_ir_{i-1}$.)

Let $\varrho:G_n\to{\Cal L}$ be a representation. Fix some $i$ and
suppose that $g:=\varrho(g_{i(i-1)})$ is hyperbolic. For~every
$t\in\Bbb R$, define a representation $\varrho E_i(t)$ as
$$\varrho E_i(t)(g_{(i+1)i}):=\varrho(g_{(i+1)i})g^{-2t},\qquad\varrho
E_i(t)(g_{(i-1)(i-2)}):=g^{2t}\varrho(g_{(i-1)(i-2)}),$$
$$\varrho E_i(t)(g_{j(j-1)}):=\varrho(g_{j(j-1)})\quad\text{for all
}j\notin\{i-1,i+1\}.$$
If $\varrho$ is induced by some $\hat\varrho:H_n\to{\Cal L}$, then
$\hat\varrho(r_i)=R(q_i)$ and $\hat\varrho(r_{i-1})=R(q_{i-1})$ for
some $q_i,q_{i-1}\in\B W$ belonging to the axis of $g$. As is easy to
see, $g^tR(q_i)g^{-t}=R(q_i)g^{-2t}$ and
$g^tR(q_{i-1})g^{-t}=g^{2t}R(q_{i-1})$. In other words, we obtain an
extension of the action of ${\Cal E}_n$ on ${\Cal R}H_n$ (and on
$\Cal H_n$) to that on ${\Cal R}G_n$ (and on $\Cal T_n$).

\bigskip

\centerline{\bf6.~References}

\medskip

[AGG]~S.~Anan$'$in, C.~H.~Grossi, N.~Gusevskii, {\it Complex Hyperbolic
Structures on Disc Bundles over Surfaces. {\rm I.}~General Settings. A
Series of Examples,} available at http://arxiv.org/abs/math/0511741

\smallskip

[AGr]~S.~Anan'in, C.~H.~Grossi, {\it Coordinate-Free Classic
Geometries. {\rm I.}~Projective Case,} available at
http://arxiv.org/abs/math/0702714

\smallskip

[Ana1]~S.~Anan$'$in, {\it Complex Hyperbolic Structures on Disc Bundles
over Surfaces. {\rm III.}~Path-following Isometry, Tectonics, and
Pentagons,} in preparation.

\smallskip

[Ana2]~S.~Anan$'$in, {\it A Hyperelliptic View on Teichm\"uller Space.
{\rm II},} in preparation.

\smallskip

[Bir]~J.~S.~Birman, {\it Braids, Links and Mapping Class Groups,}
Ann.~Math.~Studies {\bf82}, Princeton University Press, 1975.

\smallskip

[BIW]~M.~Burger, A.~Iozzi, A.~Wienhard, {\it Surface group
representations with maximal Toledo invariant,} available at
http://arxiv.org/abs/math/0605656

\smallskip

[Gol1]~W.~M.~Goldman, {\it Topological components of spaces of
representations,} Invent.~Math.~{\bf93} (1988), No.~3, 557--607.

\smallskip

[Gol2]~W.~M.~Goldman, {\it Complex Hyperbolic Geometry,} Oxford
Mathematical Monographs. Oxford Science Publications. The Clarendon
Press, Oxford University Press, New York, 1999, xx+316 pp.

\smallskip

[ImT]~Y.~Imayoshi, M.~Taniguchi, {\it An Introduction to Teichm\"uller
Spaces,} Springer-Verlag, Tokyo, 1999, x+279 pp.

\smallskip

[KMa]~V.~Koziarz, J.~Maubon, {\it Harmonic maps and representations of
non-uniform lattices of}\break

\vskip-12pt

\noindent
$\PU(m,1)$, available at http://arxiv.org/abs/math/0309193

\smallskip

[Mac]~C.~Maclachlan, {\it Smooth coverings of hyperelliptic surfaces,}
Quart.~J.~Math.~Oxford Ser.~(2) {\bf22} (1971), 117-123.

\smallskip

[Mas]~B.~Maskit, {\it Kleinian Groups,} Grundlehren der mathematischen
Wissenschaften {\bf287}, Springer-Verlag, 1987, xiv+326 pp.

\smallskip

[Stu]~M.~Stukow, {\it Small torsion generating sets for hyperelliptic
mapping class groups,} Topology and its Applications {\bf145} (2004),
83--90.

\smallskip

[Tol]~D.~Toledo, {\it Representations of surface groups in complex
hyperbolic space,} J.~Differential Geom. {\bf29} (1989), No.~1,
125--133.

\enddocument